\RequirePackage{fix-cm}
\documentclass{svjour3}                     
\smartqed  
\usepackage{graphicx}
\usepackage{amsmath}               
\usepackage{amsfonts}                 
\usepackage{varwidth}
\usepackage{parskip}
\usepackage{hyperref}
\usepackage{rotating}
\usepackage{pdflscape}
\usepackage[numbers]{natbib}
\usepackage{caption}
\usepackage{subcaption}

\newcommand{\bea}{\begin{eqnarray}}
\newcommand{\eea}{\end{eqnarray}}
\newcommand{\bee}{\begin{eqnarray*}}
\newcommand{\eee}{\end{eqnarray*}}
\newcommand{\nn}{\nonumber}
\newcommand{\lb}{\label}
\newcommand{\la}{\lambda}
\newcommand{\bpi}{\boldsymbol{\pi}}

\newcommand{\bec}{\textbf{e}}
\begin{document}

\title{Stationary queue and server content distribution of a batch-size-dependent service queue with batch Markovian arrival process: $BMAP/G^{(a,b)}_n/1$}

\author{ S. Pradhan        \and
        U.C. Gupta
}
\institute{S. Pradhan \at
               Department of Mathematics\\
Visvesvaraya National Institute of Technology, Nagpur-440010, India \\
              \email{spiitkgp11@gmail.com}           
           \and
           U.C. Gupta \at
              Department of Mathematics\\
Indian Institute of Technology, Kharagpur-721302, India\\
\email{	umesh@maths.iitkgp.ernet.in }
}

\date{Received: date / Accepted: date}

\maketitle

\begin{abstract}
Queueing systems with batch Markovian arrival process (BMAP) have paramount applications in the domain of wireless communication. The BMAP has been used to model the superposition of video sources and to approximate the super-position of data, voice and video traffic. This paper analyzes an infinite-buffer generally distributed batch-service queue with BMAP, general bulk service ($a,b$) rule and batch-size-dependent service time. In this proposed analysis, we mainly focus on deriving the bivariate vector generating function of queue and server content distribution together at departure epoch using supplementary variable technique. The mathematical procedure for the complete extraction of distribution at departure epoch has been discussed and using those extracted probabilities, we achieve the queue and server content distribution at arbitrary epoch. Finally, numerical illustrations have been carried out in order to make a deep insight to the readers which contains deterministic as well as phase-type service time distributions.
\end{abstract}
\section{Introduction}
The traditional queueing systems with standard arrival process such as Poisson, renewal or phase-type cannot adequately capture the correlation among the inter-arrivals, usually occurring in high-speed teletraffic networks such as Web browsing, VoIP, and teleconferencing. Consequently, it is necessary to have a suitable arrival process which can deal with multimedia applications over teletraffic networks. One such stochastic point process is batch Markovian arrival process (BMAP) which generalizes the standard Poisson process (and other point processes) by allowing correlated inter-arrival time of batches.\\
\hspace*{0.3cm}Lucantoni et al. \cite{lucantoni1990single} and Lucantoni \cite{lucantoni1991new} introduced BMAP which is a convenient representation of the versatile Markovian point process, see Neuts \cite{neuts1979versatile}. Later on some researchers have analyzed finite-/infinite- buffer queueing systems with Markovian arrival process (MAP) or BMAP, e.g., see Lee et al. \cite{lee2003new}, Dudin et al. \cite{dudin2005analysis}, Chaudhry et al. \cite{chaudhry2013simple}, Banik \cite{banik2015single}, Gupta et al. \cite{gupta2016alternative}. Some authors have also studied discrete Markovian/batch Markovian arrival process (D-MAP/D-BMAP), e.g., see Chaudhry et al. \cite{chaudhry2002waiting}, Samanta et al. \cite{samanta2015waiting} and references therein.\\
\hspace*{0.3cm}In recent times, a few researchers have focused on batch-service queue with batch-size-dependent service due to their wide range of applicability in production and transportation, package delivery, group testing of blood samples, telecommunication networks etc. For more detail see Bar-Lev et al. \cite{bar2007applications}, Claeys et al. \cite{claeys2010queueing,claeys2014applicability}, Pradhan et al. \cite{pradhan2015analyzing,pradhan2015queue}. In order to maximize the serving capability of the system, server content distribution plays a noteworthy role. On account of this, in a series of papers, Pradhan and Gupta \cite{gupta2015queue,pradhan2015PEVA,pradhan2017ANOR,pradhan2018IJOR} have derived bivariate probability/vector generating functions (pgf/VGF) of queue length and number with the departing batch. Furthermore, they calculated the queue and system length distributions at different epochs. On the other hand, Banerjee et al. \cite{banerjee2015analysis} discussed a finite-buffer MAP$/G^{(a,b)}_r/1/N$ queue with batch-size-dependent service and obtained joint distribution of queue length and size of the departing batch through embedded Markov chain technique (EMCT). The counter discrete version of the above model i.e., D-MAP$/G^{(1,a,b)}/1/N$ queue have been considered by Yu and Alfa \cite{yu2015algorithm} wherein they obtained joint queue length and server content distribution at various epochs by employing both EMCT and quasi-birth-and-death (QBD) process. Claeys  et al. \cite{claeys2013analysis} analyzed an infinite-buffer D-BMAP$/G^{(l,c)}_r/1$ queue with batch-size-dependent service time and focused on deriving the VGF for both queue length and server content distribution at arbitrary slot. However, they did not provide any procedure for the complete extraction of queue and server content distribution together. They also investigated the influence of correlation of the arrival process on the behavior of the system.\\
\hspace*{0.3cm}To the best of authors' knowledge, an infinite-buffer batch-service queueing model with BMAP, general bulk-service $(a,b)$ rule and batch-size-dependent service in continuous-time set-up has not yet been discussed in the literature. Moreover, the complete queue and server content distribution together as well as the \emph{only} queue-length distribution are not available in the literature so far for the concerned queue. \\
\hspace*{0.3cm}In view of this, our main objective in this paper, is to achieve the complete queue and server content distribution together for an infinite-buffer batch-service queue with BMAP and batch-size-dependent service: BMAP$/G^{(a,b)}_n/1$. The model is described in detail in the next section. Surprisingly the analysis of the concerned queue is difficult from both the modeling as well as computational point of view using EMCT mainly due to two reasons. Firstly, the construction the transition probability matrix (TPM) is quite challenging task (if not impossible). Secondly, the derivation of the analytic expression of bivariate VGF of queue and server content distribution using EMCT is tedious and difficult task, if not impossible. Fortunately, use of supplementary variable technique (SVT), where remaining service time of batch in service is taken as supplementary variable, eventually keeps us away from the use of the TPM directly. This leads the bivariate VGF of queue length and departing batch content in an effortless way, and also builds up the relationship between the probability vectors at arbitrary and departure epoch as a by-product. It may also be noted here that, the inclusion of BMAP makes the mathematical as well as computational analysis much more complex as compared to MAP because the BMAP deals with batch-arrival of customers and a selected customer from the arriving batch may be served in a different
batch from its arriving batch. Moreover, for the determination of unknown probability vectors as well as in extraction of probabilities, we deal with very complicated analytic expressions in case of BMAP. The significant contributions in this paper are: ($i$) a bivariate VGF of queue length and number in served batch is derived using SVT, ($ii$) the complete procedure of extraction of joint distribution of queue content and size of departing batch in terms of roots of the characteristic equation is provided, ($iii$) a relationship between departure and arbitrary epoch probability vectors have been generated, ($iv$) in order to manifest the feasibility and applicability of the proposed methodology and results, various numerical examples are reported with phase (PH) type as well as deterministic service time distributions.\\
\hspace*{0.3cm}The rest part of this paper is presented as follows: after describing the model in detail in the next section, we provide the governing equations of the system in Section \ref{sec3}. Section \ref{sec4} discusses the procedure of obtaining the joint distribution at departure epoch which includes derivation of bivariate VGF at departure epoch, determination of unknown vectors and finally the extraction of probability vectors. Two relations between departure and arbitrary, and arbitrary and pre-arrival epochs are obtained in Section \ref{sec5} and Section \ref{sec6}, respectively. Queue length, system length and server content distribution along with some relevant performance measures are provided in Section \ref{sec7}. Finally, Section \ref{sec8} presents numerical illustration  followed by the final conclusion.
\section{Model description and preliminaries}\lb{sec2}
\begin{itemize}
\item \textbf{Arrival process:} Customers arrive at the queueing system according to
an $m$-state batch Markovian arrival process (BMAP). Although the arrival process is well described by Lucantoni \cite{lucantoni1991new}, we give a brief description again for the sake of completeness and to clarify the notations which are used in this analysis. The arrival process is characterized by $m \times m$ matrices $\textbf{D}_k$; $k \geq 0$, where $(i,j)$-th $(1\leq i,j\leq m, i\neq j)$ element of $\textbf{D}_0$, is the state transition rate from state $i$ to state $j$ in the underlying Markov chain without an arrival, and $(i,j)$-th $(1\leq i,j\leq m, i\neq j)$ element of $\textbf{D}_k$; $k \geq 1$, is the state transition rate from state $i$ to state $j$ in the underlying Markov chain with an arrival of batch size $k$. The matrix $\textbf{D}_0$ has non-negative off-diagonal and negative diagonal elements, and the matrix $\textbf{D}_k$; $k \geq 1$, has non-negative elements. Let $N(t)$ denote the number of arrivals in $(0,t]$ and $J(t)$ be the state of the underlying Markov chain at time $t$ with state space $\{i : 1 \leq i \leq m\}$. Then $[N(t), J(t)]$ is a two-dimensional Markov process of BMAP with state space $\{(n, i) : n \geq 0, 1 \leq i \leq m\}$.\\
\hspace*{0.3cm}Further, $\textbf{D} = \sum_{k=0}^{\infty}\textbf{D}_k$ is the infinitesimal generator of the underlying Markov
chain $J(t)$. Let  $\boldsymbol{\overline \pi}=[\overline \pi_1,\overline \pi_2,\ldots,\overline \pi_m]$ be the stationary probability vector such that $ \boldsymbol{\overline \pi}\textbf{D}=\textbf{0},~~\boldsymbol{\overline \pi}\textbf{e}=1$, where \textbf{0} denotes a zero matrix of appropriate dimension and \textbf{e} is a $m \times 1$ column vector with all elements as 1. Throughout the analysis we use these notations, but when needed, the dimension of \textbf{0} and \textbf{e} will be identified with a suffix. The fundamental arrival rate of the above Markov process is given by $\la^*=\boldsymbol{\overline \pi}\textbf{D}\textbf{e}$. We assume that the notation $\textbf{I}$ stands for an identity matrix of appropriate dimension.\\
\hspace*{0.3cm}Then the average arrival rate $\la^*$ and average batch arrival rate $\la_g$ of the stationary BMAP are given by $ \la^*=\overline\bpi \sum_{k=1}^{\infty}k \textbf{D}_k \textbf{e}$,
$\la_g=\overline\bpi  \sum_{k=1}^{\infty} \textbf{D}_k \textbf{e}= \overline\bpi \textbf{ D}'_1$,
respectively, where $\textbf{D}'_n =  \sum_{k=n}^{\infty}\textbf{D}_k,~~n \geq1$.

\item \textbf{Service discipline:} The customers are served in group/batches according to general bulk service ($a,b$) rule. The server only starts service if the queue contains at least as many customers as the service threshold $`a$'. If queue contains fewer than $`a$' customers, then the server is said to be in idle period. For the queue size $r$ ($a\leq r\leq b$), entire group of customers are taken for service. When the queue size exceeds $`b$', then the server can process maximum $`b$' customers and others remain in the queue for the next round of service.

\item \textbf{Service process:} A service time is the length of a service period, and the consecutive service times are independently and identically distributed. We assume that the service times of the batches follow general distribution and to be dependent on batch size of ongoing service. The length of random service time of a batch of size $r$ is denoted by the random variable $T_r$ ($a \leq r \leq b$) with probability density function (pdf) $s_r(t)$, distribution function $S_r(t)$, the Laplace-Stieltjes transform (L.-S.T.) $\widetilde{S}_r(\theta)$ and the mean service time $ \frac{1}{\mu_r}=s_r=-\widetilde{S}_r^{(1)}(0)$, where $\widetilde{S}_r^{(1)}(0)$ is the derivative of $\widetilde{S}_r(\theta)$ evaluated at $\theta=0$.

\item The traffic intensity of the system is given by $\rho= \frac{\la^*}{b\mu_b} $ and $\rho<1$ ensures the stability of the system.

\end{itemize}
\section{Governing equations of the system}\lb{sec3}
In this section, we develop the governing equations of the model where SVT is employed with remaining service time of a batch as supplementary variable. We define the state of the system at time $t$ as
\begin{itemize}

\item $N_q(t)$ $\equiv$ Number of customers in the queue waiting for service,

\item $S(t)$ $\equiv$ Number of customers with the server,

\item $J(t)$ $\equiv$ State of the underlying chain of the $MAP$ and

\item $U(t)$ $\equiv$ Remaining service time of a batch in service (if any).

\end{itemize}
Let us define for $1\leq i \leq m$,
\begin{small}
\bea p_i(n,0;t) &=& \mbox{Pr}\{ N_q(t)=n, S(t)=0, J(t)=i, \mbox{server idle}\},~~0\leq n \leq a-1, \nn\\
\pi_i(n,r,u;t)du &=& \mbox{Pr}\{ N_q(t)=n, S(t)=r, J(t)=i, u<U(t)\leq u+du, \mbox{server busy} \},\nn\\
&&\hspace{5cm}n\geq 0,~ a\leq r \leq b, ~u\geq 0.\nn\eea
\end{small}
Also let us define the limiting probabilities as
\bea p_i(n,0) &=& \lim_{t\rightarrow \infty} p_i(n,0;t), ~~ 1\leq i \leq m,~0\leq n \leq a-1, \nn\\
\pi_i(n,r,u) &=&  \lim_{t\rightarrow \infty} \pi_i(n,r,u;t),~~ 1\leq i \leq m, ~n\geq 0,~ a\leq r \leq b. \nn \eea
Let us define the probability vectors $\textbf{p}(n,0) = \left(p_1(n,0),\ldots,p_m(n,0)\right)$ and\\ $\boldsymbol{\pi}(n,r,u) = \left(\pi_1(n,r,u),\ldots,\pi_m(n,r,u)\right)$.
We relate the states of the system at two consecutive times $t$ and $t+dt$, and by considering each phase, in steady-state, we obtain the following equations in vector and matrix form:
\begin{small}
\bea 0 &=& \textbf{p}(0,0)\boldsymbol{D}_0+\sum_{j=a}^{b}\boldsymbol{\pi}(0,j,0)\lb{map7}\\
0 &=& \textbf{p}(n,0)\boldsymbol{D}_0+\sum_{i=1}^{n}\textbf{p}(n-i,0)\boldsymbol{D}_i+\sum_{j=a}^{b}\boldsymbol{\pi}(n,j,0),~1\leq n \leq a-1\lb{map8}\\
-\frac{d}{du}\boldsymbol{\pi}(0,r,u)&=&\boldsymbol{\pi}(0,r,u)\boldsymbol{D}_0+\sum_{i=0}^{a-1}\textbf{p}(i,0)\boldsymbol{D}_{r-i}s_r(u)
+\sum_{j=a}^{b}\boldsymbol{\pi}(r,j,0)s_r(u),~~a\leq r \leq b \lb{map9}\\
-\frac{d}{du}\boldsymbol{\pi}(n,r,u) &=& \boldsymbol{\pi}(n,r,u)\boldsymbol{D}_0+\sum_{i=1}^{n}\boldsymbol{\pi}(n-i,r,u)\boldsymbol{D}_i,~ a\leq r \leq b-1,~n\geq 1\lb{map11}\\
-\frac{d}{du}\boldsymbol{\pi}(n,b,u)&=&\boldsymbol{\pi}(n,b,u)\boldsymbol{D}_0+\sum_{i=1}^{n}\boldsymbol{\pi}(n-i,b,u)\boldsymbol{D}_i
+\sum_{j=a}^{b}\textbf{p}(b-j,0)\boldsymbol{D}_{n+j}s_b(u)\nn\\
&&~~~~~~~~~~~~~~~~~~~~~~~~~~~~~~~~~~~~~+\sum_{j=a}^{b}\boldsymbol{\pi}(n+b,j,0)s_b(u),~ n\geq 1\lb{map12}.\eea
\end{small}
~\\
\noindent Our main objective is to achieve the joint distribution of queue content and number with the departing batch and phase of the arrival process from the entire set of governing difference-differential equations (\ref{map7}) to (\ref{map12}) of the model under consideration. In view of this, let us define the Laplace transform of $\boldsymbol{\pi}(n,r,u)$ as
\bea \widetilde{\boldsymbol{\pi}}(n,r,\theta)&=&\int_{0}^{\infty}e^{-\theta u}\boldsymbol{\pi}(n,r,u)du,\quad a\leq r \leq b,~n\geq 0,~ \mathfrak{R}(\theta)\geq 0 \\
\mbox{so that}, ~~\boldsymbol{\pi}(n,r)&=& \widetilde{\boldsymbol{\pi}}(n,r,0)= \int_{0}^{\infty}\boldsymbol{\pi}(n,r,u)du,~ a\leq r \leq b,~n\geq 0\lb{map13}.\eea
Now the equations (\ref{map9}) to (\ref{map12}) are transformed through the multiplication by $e^{-\theta u}$ and integrating with respect to $u$ over 0 to $\infty$ and are given by:
\begin{small}
\bea \hspace*{-1.0cm}-\theta \widetilde{\boldsymbol{\pi}}(0,r,\theta)+\boldsymbol{\pi}(0,r,0)&=&\widetilde{\boldsymbol{\pi}}(0,r,\theta)\textbf{D}_0
+\sum_{i=0}^{a-1}\textbf{p}(i,0)\textbf{D}_{r-i}\widetilde{S}_r(\theta)+\sum_{j=a}^{b}\boldsymbol{\pi}(r,j,0)\widetilde{S}_r(\theta),~a\leq r \leq b \lb{map14}\\
\hspace*{-1.0cm}-\theta \widetilde{\boldsymbol{\pi}}(n,r,\theta)+\boldsymbol{\pi}(n,r,0) &=& \widetilde{\boldsymbol{\pi}}(n,r,\theta)\textbf{D}_0
+\sum_{i=1}^{n}\widetilde{\boldsymbol{\pi}}(n-i,r,\theta)\textbf{D}_i,~ a\leq r \leq b-1,~n\geq 1\lb{map16}\\
\hspace*{-1.0cm}-\theta \widetilde{\boldsymbol{\pi}}(n,b,\theta)+\boldsymbol{\pi}(n,b,0)&=&\widetilde{\boldsymbol{\pi}}(n,b,\theta)\textbf{D}_0
+\sum_{i=1}^{n}\widetilde{\boldsymbol{\pi}}(n-i,b,\theta)\textbf{D}_{i}+\sum_{j=a}^{b}\textbf{p}(b-j,0)\textbf{D}_{n+j}\widetilde{S}_b(\theta)\nn\\
&&\hspace*{3cm}+\sum_{j=a}^{b}\boldsymbol{\pi}(n+b,j,0) \widetilde{S}_b(\theta),~n\geq 1\lb{map17} \eea
\end{small}
Our center of focus is to acquire the joint distribution of queue content as well as server
content at departure and arbitrary epoch. For this purpose, first we define the following joint probabilities at
departure epoch when arrival process is in phase $i$ as
\bea \pi_i^+(n,r)&\equiv& \mbox{joint probability that there are}~ n ~(n\geq 0)~ \mbox{customers in the queue} \nn\\
&& \mbox{and the arrival process is in phase}~ i ~(1\leq i\leq m)~ \mbox{immediately after} \nn\\
&& \mbox{the departure of a batch of size}~ r~ (a \leq r \leq b),\lb{map18}\\
\psi_i^+(n)&\equiv&\mbox{Pr}\{\mbox{queue contains}~ n ~\mbox{customers}~ \mbox{and the arrival process is in phase}~ i~\nn\\
&&~~~~~~ \mbox{at departure epoch of a batch}\}\nn\\
&=&\sum_{r=a}^{b}\pi_i^+(n,r),\lb{map19}\\
\phi_i^+(r)&\equiv&\mbox{Pr}\{\mbox{there are}~ r ~\mbox{customers with the departing batch}~ \mbox{and the arrival }\nn\\
&&\mbox{~~~~process is in phase}~ i~\}~\nn \\
&=&\sum_{n=0}^{\infty}\pi_i^+(n,r).\lb{map20}\eea
As a consequence, we have the joint probability vectors as
\bea \boldsymbol{\pi}^+(n,r)&=&\left[\pi_1^+(n,r),\ldots,\pi_m^+(n,r)\right]\lb{map21}\\
\boldsymbol{\psi}^+(n)&=&\left[\psi_1^+(n),\ldots,\psi_m^+(n)\right]\lb{map22}\\
\boldsymbol{\phi}^+(r)&=&\left[\phi_1^+(r),\ldots,\phi_m^+(r)\right]\lb{map23}\eea
Now we propose some results which will be used in the analysis.
\begin{lemma}
The probability vectors $\boldsymbol{\pi}^+(n,r)$ and $\boldsymbol{\pi}(n,r,0)$ are connected by the relation
\bea \boldsymbol{\pi}^+(n,r)=\frac{\boldsymbol{\pi}(n,r,0)}{\displaystyle\sum_{j=0}^{\infty}\sum_{l=a}^{b}\boldsymbol{\pi}(j,l,0)\bec}\lb{map24}\eea
\end{lemma}\lb{maplm1}
\noindent\emph{Proof:} As $\boldsymbol{\pi}^+(n,r)$ and $\boldsymbol{\pi}(n,r,0)$ differ by a constant, using  $\displaystyle\sum_{n=0}^{\infty}\sum_{r=a}^{b}\boldsymbol{\pi}^+(n,r)\bec = 1$, we are led to the desired result.
~\\
\begin{lemma}
The probability vectors $\textbf{p}(n,0)$ and $\bpi(n,r,0)$ are related by
\begin{small}
\bea \sum_{n=0}^{\infty}\sum_{r=a}^{b}\bpi(n,r,0)\bec=\frac{1-\displaystyle\sum_{n=0}^{a-1}\textbf{p}(n,0)\bec}{\omega} \lb{bmap19}\eea
where
\begin{small}
\bea \omega&=&\sum_{n=0}^{a-1}\boldsymbol{\psi}^+(n) \left(\sum_{\ell=a}^{b}\boldsymbol{C}_{\ell n}s_{\ell}+\sum_{j=n}^{a-1}\boldsymbol{M}_{j n}\sum_{i=b+1-j}^{\infty}\boldsymbol{\overline{D}}_i s_{b}\right) \bec+\sum_{n=a}^{b}s_n \boldsymbol{\psi}^+(n)\bec\nn\\
&&\hspace*{6.0cm}+s_b\sum_{n=b+1}^{\infty}\boldsymbol{\psi}^+(n)\bec \lb{bmap2}\\
\boldsymbol{\overline{D}}_i&=&\left(-\boldsymbol{D}_0\right)^{-1}\boldsymbol{D}_i\nn\eea
\end{small}
\mbox{and}\\
\\
$\boldsymbol{C}_{\ell,n} = \left\{\begin{array}{r@{\mskip\thickmuskip}l}
& \boldsymbol{\overline{D}}_{\ell-n},~~~~~~~~~~~~~~~~~~~~~~~~~~~~n=a-1,\\
&~~\\
&\displaystyle\sum_{j=n+1}^{a-1}\boldsymbol{\overline{D}}_{j-n}\boldsymbol{C}_{\ell,j}+\boldsymbol{\overline{D}}_{\ell-n},~~n=0,1,\ldots,a-2,~\ell=a,\ldots,b.
 \end{array}\right.$\\
 ~~\\
\mbox{and}\\
 $\boldsymbol{M}_{n,i}=\displaystyle\sum_{j=i+1}^{n-1}M_{n,j}\boldsymbol{\overline{D}}_{j-i}+\boldsymbol{\overline{D}}_{n-i},~i=0,1,\ldots,n-2$, \mbox{with} $\boldsymbol{M}_{n,n-1}=\boldsymbol{\overline{D}}_1$ \mbox{and} $\boldsymbol{M}_{n,n}=\textbf{I}$

\end{small}
\end{lemma}\lb{maplm2}
\noindent\emph{Proof:} Post multiplying (\ref{map14}) - (\ref{map17}) by the vector $\textbf{e}$, adding them and using $\textbf{D}\textbf{e}$=\textbf{0}, we obtain
\begin{small}
\bea  \sum_{n=0}^{\infty}\sum_{r=a}^{b}\widetilde{\boldsymbol{\pi}}(n,r,\theta)\textbf{e}&=&
\sum_{n=0}^{a-1}\sum_{r=a}^{b}\bpi(n,r,0)\left[\frac{\textbf{I}-\left(\displaystyle\sum_{\ell=a}^{b}\boldsymbol{C}_{\ell n}\widetilde{S}_{\ell}(\theta)+\sum_{j=n}^{a-1}\boldsymbol{M}_{j n}\sum_{i=b+1-j}^{\infty}\boldsymbol{\overline{D}}_i\widetilde{S}_{b}(\theta)\right)}{\theta}\right]\bec\nn\\
&&+\sum_{n=a}^{b}\sum_{r=a}^{b}\frac{1-\widetilde{S}_n(\theta)}{\theta}\bpi(n,r,0)\bec
+\frac{1-\widetilde{S}_b(\theta)}{\theta}\sum_{n=b+1}^{\infty}\sum_{r=a}^{b}\bpi(n,r,0)\bec \nn\eea
\end{small}
Taking limit $\theta \rightarrow 0$ in the above expression, using l'H$\hat{o}$spital's rule, equation (\ref{map24}) and the normalizing condition
\bea \sum_{n=0}^{a-1}\textbf{p}(n,0)\bec + \sum_{n=0}^{\infty}\sum_{r=a}^{b}\bpi(n,r)\bec &=& 1 \lb{map27}\eea
we obtain the desired result as (\ref{bmap19}).
~\\
\begin{lemma}
The value of $\displaystyle \sum_{n=0}^{\infty}\sum_{r=a}^{b}\bpi(n,r,0)\bec$ is given as
\begin{small}
\bea \sum_{n=0}^{\infty}\sum_{r=a}^{b}\bpi(n,r,0)\bec =\frac{1-\displaystyle \sum_{n=0}^{a-1}\textbf{p}(n,0)\bec}{\omega}
=\frac{1}{\omega+\displaystyle \sum_{n=0}^{a-1}\sum_{j=0}^{n}\boldsymbol{\psi}^+(j)\boldsymbol{M}_{n,j}(-\textbf{D}_0)^{-1}\bec}\lb{map28}\eea
\end{small}
where $\omega$ is presented in equation (\ref{bmap2}).
\end{lemma}\lb{maplm3}
\noindent\emph{Proof}: Dividing (\ref{map7}) by $\displaystyle \sum_{n=0}^{\infty}\sum_{r=a}^{b}\bpi(n,r,0)\bec$ and using (\ref{map24}), (\ref{map22}) and (\ref{bmap19}) we obtain
\bea \textbf{p}(0,0)= \frac{1-\displaystyle \sum_{n=0}^{a-1}\textbf{p}(n,0)\bec}{\omega}\boldsymbol{\psi}^+(0)(-\textbf{D}_0)^{-1}\lb{map30}\eea
Similarly, from (\ref{map8}), and using (\ref{bmap19}) we obtain
\bea \textbf{p}(n,0)= \frac{1-\displaystyle\sum_{n=0}^{a-1}\textbf{p}(n,0)\bec}{\omega}\displaystyle\sum_{j=0}^{n}\boldsymbol{\psi}^+(j)\boldsymbol{M}_{n,j}(-\textbf{D}_0)^{-1},~0\leq n\leq a-1 \lb{map31}\eea
Post multiplying (\ref{map31}) by $\bec$ we obtain
\bea \textbf{p}(n,0)\bec= \frac{1-\displaystyle \sum_{n=0}^{a-1}\textbf{p}(n,0)\bec}{\omega}\displaystyle\sum_{j=0}^{n}\boldsymbol{M}_{n,j}(-\textbf{D}_0)^{-1}\bec, ~~0\leq n\leq a-1 \nn\eea
which implies that
\bea \frac{1-\displaystyle\sum_{n=0}^{a-1}\textbf{p}(n,0)\bec}{\omega}
=\displaystyle \frac{1}{\omega+\displaystyle\sum_{n=0}^{a-1}\sum_{j=0}^{n}\boldsymbol{\psi}^+(j)\boldsymbol{M}_{n,j}(-\textbf{D}_0)^{-1}\bec}. \lb{map32}\eea
Combining (\ref{bmap19}) and (\ref{map32}) we obtain the desired result of lemma \ref{maplm3}.
\section{Distribution of queue length and number with the departing batch}\lb{sec4}
The purpose of this section is three fold: (i) the derivation of bivariate VGF of queue-length and number with departing batch, (ii) the determination of unknown probability vectors appearing in the numerator of VGF, (iii) the procedure of extraction of the complete joint distribution of queue and server content at departure epoch. On account of this, we first proceed with the derivation of bivariate VGF of queue length and departing batch content in the following subsection.
\subsection{Bivariate VGF at departure epoch}
First we define the pgfs of sequences $\{\widetilde{\pi}_i(n,r,\theta)\}$, $\{\pi^+_i(n,r)\}$ and $\{\psi_i^+(n)\}$ when arrival process is in phase $i$:
\bea \Pi_i(z,y,\theta) &\equiv& \sum_{n=0}^{\infty}\sum_{r=a}^{b}\widetilde{\pi}_i(n,r,\theta)z^n y^r ,~~1\leq i \leq m,~ |z|\leq 1,~|y|\leq 1, \lb{map33}\\
\Pi^+_i(z,y)&\equiv& \sum_{n=0}^{\infty}\sum_{r=a}^{b}\pi^+_i(n,r)z^n y^r ,~~1\leq i \leq m,~ |z|\leq 1,~|y|\leq 1, \lb{map34}\\
\Psi_i^+(z)&\equiv&\sum_{n=0}^{\infty}\psi_i^+(n)z^n\equiv \Pi^+_i(z,1),~~1\leq i \leq m,~ |z|\leq 1,~|y|\leq 1. \lb{map35}\eea
Consequently we have the following VGF as:
\bea \boldsymbol{\Pi}(z,y,\theta)&=&\left[\Pi_1(z,y,\theta),\ldots,\Pi_m(z,y,\theta)\right]\lb{map36}\\
\boldsymbol{\Pi}^+(z,y)&=&\left[\Pi^+_1(z,y),\ldots,\Pi^+_m(z,y)\right]\lb{map37}\\
\boldsymbol{\Psi}^+(z)&=&\left[\Psi^+_1(z),\ldots,\Psi^+_m(z)\right] \lb{map38}.\eea
Now, the transformed equations (\ref{map14}) - (\ref{map17}) are multiplied by appropriate power of $z$ and $y$ in order to derive the VGF. Then summing over $n$ from 0 to $\infty$, and $r$ from $a$ to $b$, and using (\ref{map33}) and (\ref{map36}) we get
\begin{small}
\bea \hspace*{-0.5cm}\boldsymbol{\Pi}(z,y,\theta)\left[-\theta \textbf{I}-\textbf{D}(z)\right]&=& \sum_{n=0}^{a-1}\sum_{j=a}^{b}\left\{\sum_{\ell=0}^{n}\sum_{r=a}^{b}\bpi(\ell,r,0)\boldsymbol{M}_{n,\ell}\right\}\overline{\textbf{D}}_{j-n}~y^j
\widetilde{S}_j(\theta)\nn\\
&&+\sum_{n=a}^{b}\sum_{r=a}^{b}\bpi(n,r,0)y^n\widetilde{S}_n(\theta)\nn\\
&&+\sum_{n=0}^{a-1}\sum_{\ell=0}^{n}\sum_{r=a}^{b}\bpi(\ell,r,0)\boldsymbol{M}_{n,\ell}\sum_{i=b+1-n}^{\infty}\overline{\textbf{D}}_i~ z^{i-b+n} y^b \widetilde{S}_b(\theta)\nn\\
&&+\sum_{n=b+1}^{\infty}\sum_{r=a}^{b}\bpi(n,r,0)z^{n-b} y^b\widetilde{S}_b(\theta)-\sum_{n=0}^{\infty}\sum_{r=a}^{b}\bpi(n,r,0)z^n y^r \lb{map39}\eea
\end{small}
From the above equation, our principal aim is to achieve the bivariate VGF of queue length and size of the departing batch which can be accomplished by making left hand side of (\ref{map39}) to zero. In order to accomplish this, the eigenvalues and eigenvectors of $-\textbf{D}(z)$ have to be used, see Lee et al. \cite{lee2001decompositions}. Let us denote its eigenvalues as $\alpha_1(z),\ldots,\alpha_m(z)$, and $\xi_1(z),\ldots,\xi_m(z)$ as the right eigenvectors. Therefore we have
\bea -\textbf{D}(z)\xi_i(z)=\alpha_i(z)\xi_i(z) \nn\eea which implies that
\bea \left[\alpha_i(z)\textbf{I}+\textbf{D}(z)\right]\xi_i(z)=0\lb{map40}\eea
Now substituting $\theta=\alpha_i(z)$ in (\ref{map39}) and post multiplying both sides by $\xi_i(z)$, we obtain
\begin{small}
\bea \hspace*{-0.5cm}\boldsymbol{\Pi}(z,y,\alpha_i(z))\left[-\alpha_i(z)\textbf{I}-\textbf{D}(z)\right]\xi_i(z)&=&\left[ \sum_{n=0}^{a-1}\sum_{j=a}^{b}\left\{\sum_{\ell=0}^{n}\sum_{r=a}^{b}\bpi(\ell,r,0)\boldsymbol{M}_{n,\ell}\right\}\overline{\textbf{D}}_{j-n}~y^j
\widetilde{S}_j(\alpha_i(z))\right.\nn\\
&&\left.+\sum_{n=a}^{b}\sum_{r=a}^{b}\bpi(n,r,0)y^n\widetilde{S}_n(\alpha_i(z))\right.\nn\\
&&\left.+\sum_{n=0}^{a-1}\sum_{\ell=0}^{n}\sum_{r=a}^{b}\bpi(\ell,r,0)\boldsymbol{M}_{n,\ell}\sum_{j=b+1-n}^{\infty}\overline{\textbf{D}}_j~ z^{j-b+n} y^b \widetilde{S}_b(\alpha_i(z))\right.\nn\\
&&\left.+\sum_{n=b+1}^{\infty}\sum_{r=a}^{b}\bpi(n,r,0)z^{n-b} y^b\widetilde{S}_b(\alpha_i(z))\right.\nn\\
&&\left.-\sum_{n=0}^{\infty}\sum_{r=a}^{b}\bpi(n,r,0)z^n y^r \right]\xi_i(z)\lb{map41}\eea
\end{small}
Using (\ref{map40}) in (\ref{map41}) we get
\begin{small}
\bea \sum_{n=0}^{\infty}\sum_{r=a}^{b}\bpi(n,r,0)z^n y^r\xi_i(z)&=& \sum_{n=0}^{a-1}\sum_{j=a}^{b}\left\{\sum_{\ell=0}^{n}\sum_{r=a}^{b}\bpi(\ell,r,0)\boldsymbol{M}_{n,\ell}\right\}\overline{\textbf{D}}_{j-n}~y^j
\widetilde{S}_j(\alpha_i(z))\xi_i(z)\nn\\
&&+\sum_{n=a}^{b}\sum_{r=a}^{b}\bpi(n,r,0)y^n\widetilde{S}_n(\alpha_i(z))\xi_i(z)\nn\\
&&+\sum_{n=0}^{a-1}\sum_{\ell=0}^{n}\sum_{r=a}^{b}\bpi(\ell,r,0)\boldsymbol{M}_{n,\ell}\sum_{j=b+1-n}^{\infty}\overline{\textbf{D}}_j~ z^{j-b+n} y^b \widetilde{S}_b(\alpha_i(z))\xi_i(z)\nn\\
&&+\sum_{n=b+1}^{\infty}\sum_{r=a}^{b}\bpi(n,r,0)z^{n-b} y^b\widetilde{S}_b(\alpha_i(z)) \xi_i(z) \lb{map42}\eea
\end{small}
Since the equation (\ref{map42}) is true for all eigenvalues $\alpha_1(z),\ldots,\alpha_m(z)$. Consequently we have
\begin{small}
\bea \sum_{n=0}^{\infty}\sum_{r=a}^{b}\bpi(n,r,0)z^n y^r\left[\xi_1(z),\ldots,\xi_m(z)\right]&=&\sum_{n=0}^{a-1}\sum_{j=a}^{b}\left\{\sum_{\ell=0}^{n}\sum_{r=a}^{b}\bpi(\ell,r,0)\boldsymbol{M}_{n,\ell}\right\}
\overline{\textbf{D}}_{j-n}~y^j\nn\\
&&\times\left[\widetilde{S}_j(\alpha_1(z))\xi_1(z),\ldots,\widetilde{S}_j(\alpha_m(z))\xi_m(z)\right]\nn\\
&&+\sum_{n=a}^{b}\sum_{r=a}^{b}\bpi(n,r,0)y^n \nn\\
&&\times \left[\widetilde{S}_n(\alpha_1(z))\xi_1(z),\ldots,\widetilde{S}_n(\alpha_m(z))\xi_m(z)\right]\nn\\
&&+\sum_{n=0}^{a-1}\sum_{\ell=0}^{n}\sum_{r=a}^{b}\bpi(\ell,r,0)\boldsymbol{M}_{n,\ell}\sum_{j=b+1-n}^{\infty}\overline{\textbf{D}}_j~ z^{j-b+n} y^b \nn\\ &&\times\left[\widetilde{S}_b(\alpha_1(z))\xi_1(z),\ldots,\widetilde{S}_b(\alpha_m(z))\xi_m(z)\right]\nn\\
&&+\sum_{n=b+1}^{\infty}\sum_{r=a}^{b}\bpi(n,r,0)z^{n-b} y^b \nn\\
&&\times \left[\widetilde{S}_b(\alpha_1(z))\xi_1(z),\ldots,\widetilde{S}_b(\alpha_m(z))\xi_m(z)\right] \lb{map43}\eea
\end{small}
Let us define an ($m\times m$) matrix $\Delta(z)$ as follows:
\bea \Delta(z)=\left(\xi_1(z),\ldots,\xi_m(z)\right) \lb{map44}\eea
The inverse of $\Delta(z)$ exists whenever each eigenvalue is of multiplicity 1, see Nishimura \cite{nishimura1998eigenvalue}.
Then from (\ref{map43}) it follows that
\begin{small}
\bea \hspace*{-1.0cm}\sum_{n=0}^{\infty}\sum_{r=a}^{b}\bpi(n,r,0)z^n y^r&=& \sum_{n=0}^{a-1}\sum_{j=a}^{b}\left\{\sum_{\ell=0}^{n}\sum_{r=a}^{b}\bpi(\ell,r,0)\boldsymbol{M}_{n,\ell}\right\}
\overline{\textbf{D}}_{j-n}~y^j \nn\\
&&\times \Delta(z)\left[diag\left\{\widetilde{S}_j(\alpha_i(z))\right\}_{i=1}^{m}\right]\Delta^{-1}(z)\nn\\
&&+\sum_{n=a}^{b}\sum_{r=a}^{b}\bpi(n,r,0)y^n \times \Delta(z)\left[diag\left\{\widetilde{S}_n(\alpha_i(z))\right\}_{i=1}^{m}\right]\Delta^{-1}(z) \nn\\
&&+ \sum_{n=0}^{a-1}\sum_{\ell=0}^{n}\sum_{r=a}^{b}\bpi(\ell,r,0)\boldsymbol{M}_{n,\ell}\sum_{j=b+1-n}^{\infty}\overline{\textbf{D}}_j~ z^{j-b+n} y^b \nn\\
&&\times \Delta(z)\left[diag\left\{\widetilde{S}_b(\alpha_i(z))\right\}_{i=1}^{m}\right]\Delta^{-1}(z)\nn\\
&&+\sum_{n=b+1}^{\infty}\sum_{r=a}^{b}\bpi(n,r,0)z^{n-b} y^b
\times \Delta(z)\left[diag\left\{\widetilde{S}_b(\alpha_i(z))\right\}_{i=1}^{m}\right]\Delta^{-1}(z) \lb{map45}\eea
\end{small}
where $diag\left\{\widetilde{S}_r(\alpha_i(z))\right\}_{i=1}^{m}$, ($a \leq r \leq b$), is a diagonal matrix of order `$m$' with diagonal entries\\ $\widetilde{S}_r(\alpha_1(z)),\ldots,\widetilde{S}_r(\alpha_m(z))$, i.e., \\
$diag\left\{\widetilde{S}_r(\alpha_i(z))\right\}_{i=1}^{m}=\bordermatrix{ &           &             &        &    \cr
                                                                           & \widetilde{S}_r(\alpha_1(z))&  0 & \ldots &0 \cr
                                                                           & 0&  \widetilde{S}_r(\alpha_2(z)) & \ldots     &  0  \cr
                                                                           & \vdots    & \vdots      & \ddots & \vdots   \cr
                                                                           & 0&  0 & \ldots & \widetilde{S}_r(\alpha_m(z))   \cr}$
\\
\\
\noindent Now dividing (\ref{map45}) by $\displaystyle \sum_{\ell=0}^{\infty}\sum_{j=a}^{b}\boldsymbol{\pi}(\ell,j,0)\bec$ and using (\ref{map24}), (\ref{map22}), (\ref{map34}) and (\ref{map37}), we get


\begin{small}
\bea \boldsymbol{\Pi}^+(z,y)&=& \sum_{n=0}^{a-1}\sum_{j=a}^{b}\left\{\sum_{\ell=0}^{n}\boldsymbol{\psi}^+{(\ell)}\boldsymbol{M}_{n,\ell}\right\}\overline{\textbf{D}}_{j-n}~
y^j \times\Delta(z)\left[diag\left\{\widetilde{S}_j(\alpha_i(z))\right\}_{i=1}^{m}\right]\Delta^{-1}(z) \nn\\
&&+ \sum_{n=a}^{b}\boldsymbol{\psi}^+{(n)} y^n \times\Delta(z)\left[diag\left\{\widetilde{S}_n(\alpha_i(z))\right\}_{i=1}^{m}\right]\Delta^{-1}(z)  \nn\\
&&+ \sum_{n=0}^{a-1}\sum_{\ell=0}^{n}\boldsymbol{\psi}^+{(\ell)}\boldsymbol{M}_{n,\ell}\sum_{i=b+1-n}^{\infty}\overline{\textbf{D}}_i~ z^{i-b+n} y^b \times\Delta(z)\left[diag\left\{\widetilde{S}_b(\alpha_i(z))\right\}_{i=1}^{m}\right]\Delta^{-1}(z)\nn\\
&&+ \sum_{n=b+1}^{\infty}\boldsymbol{\psi}^+{(n)} z^{n-b} y^b \times\Delta(z)\left[diag\left\{\widetilde{S}_b(\alpha_i(z))\right\}_{i=1}^{m}\right]\Delta^{-1}(z) \lb{map46}\eea
\end{small}
Now using the known results of matrix algebra, the matrix $-\textbf{D}(z)$ can be written in terms of it's eigenvalues and eigenvectors as
\bea -\textbf{D}(z)&=&\Delta(z)\left[diag\left\{\alpha_i(z)\right\}_{i=1}^{m}\right]\Delta^{-1}(z) \nn\\
\left[-\textbf{D}(z)\right]^{-1}&=&\Delta(z)\left[diag\left\{\frac{1}{\alpha_i(z)}\right\}_{i=1}^{m}\right]\Delta^{-1}(z)\nn\eea
Let us denote the matrix $A^{(r)}_k=[A^{(r)}_k(x)]_{ij},~k\geq 0,~a\leq r \leq b,~1\leq i,j\leq m,~x\geq 0$, to be the conditional probability that, a departure which left at least $a$ customers in the queue with the arrival process in state $i$, the next departure occurs no later than time $x$, and during the service period of $r$ customers exactly $k$ new customers arrive; the phase of the arrival process is in phase $j$ at the departure epoch. Let us denote $\textbf{A}^{(r)}(z)$ to be the corresponding matrix generating function of the matrix $A^{(r)}_k$. Therefore we have
\bea \textbf{A}^{(r)}(z)=\sum_{j=0}^{\infty}A^{(r)}_k z^k=\int_{0}^{\infty}e^{\textbf{D}(z)t}dS_r(t) \nn\eea
Further, from matrix exponential theory we have
\bea e^{\textbf{D}(z)t}=\sum_{n=0}^{\infty}\frac{\left[\textbf{D}(z)t\right]^n}{n!}
=\Delta(z)\left[diag\left\{e^{-\alpha_i(z)t}\right\}_{i=1}^{m}\right]\Delta^{-1}(z)\nn\eea
which implies that
\bea  \textbf{A}^{(r)}(z) =\Delta(z)\left[diag\left\{\widetilde{S}_r(\alpha_i(z))\right\}_{i=1}^{m}\right]\Delta^{-1}(z),~~a\leq r \leq b\nn\eea
Therefore, from (\ref{map46}) we get

\begin{small}
\bea \boldsymbol{\Pi}^+(z,y)&=& \sum_{n=0}^{a-1}\sum_{j=a}^{b}\left\{\sum_{\ell=0}^{n}\boldsymbol{\psi}^+{(\ell)}\boldsymbol{M}_{n,\ell}\right\}\overline{\textbf{D}}_{j-n}
y^j \textbf{A}^{(j)}(z) + \sum_{n=a}^{b}\boldsymbol{\psi}^+{(n)} y^n \textbf{A}^{(n)}(z)  \nn\\&&+ \sum_{n=0}^{a-1}\sum_{\ell=0}^{n}\boldsymbol{\psi}^+{(\ell)}\boldsymbol{M}_{n,\ell}\sum_{i=b+1-n}^{\infty}\overline{\textbf{D}}_i~ z^{i-b+n} y^b \textbf{A}^{(b)}(z)\nn\\
&&+ \sum_{n=b+1}^{\infty}\boldsymbol{\psi}^+{(n)} z^{n-b} y^b \textbf{A}^{(b)}(z) \lb{map47}\eea
\end{small}
Now substituting $y=1$ in (\ref{map47}), using (\ref{map38}) and after little bit simplification, we obtain
\begin{small}
\bea \hspace*{-0.4cm}\boldsymbol{\Psi}^+(z)&=& \left[\sum_{n=0}^{a-1}\boldsymbol{\psi}^+(n)\left\{ z^b \sum_{j=a}^{b} \boldsymbol{C}_{j n} \textbf{A}^{(j)}(z)-\left(z^n-\sum_{i=b-a+2}^{b+1-n}\tau_i(z)~\boldsymbol{Q}_{b+1-n-i}~z^{b+1-i}\right)\textbf{A}^{(b)}(z)\right\}\right.\nn\\ &&\left.+\sum_{n=a}^{b-1}\boldsymbol{\psi}^+(n)\left(z^b\textbf{A}^{(n)}(z)-z^n\textbf{A}^{(b)}(z)\right)\right]
\left[z^b\textbf{I}-\textbf{A}^{(b)}(z)\right]^{-1} \lb{map48}\eea
\end{small}
where $\boldsymbol{Q}_{n} = \sum_{j=1}^{n-1}\boldsymbol{Q}_{j}\boldsymbol{\overline{D}}_{n-j}+\boldsymbol{\overline{D}}_{n},~~n=2,3,\ldots,(a-1), \mbox{with}~~\boldsymbol{Q}_0=\textbf{I}, ~\boldsymbol{Q}_1=\boldsymbol{\overline{D}}_1$ and $\tau_j(z)=\sum_{k=j}^{\infty}\boldsymbol{\overline{D}}_k z^k.$\\
Making use of (\ref{map48}) in (\ref{map47}) and use of (\ref{map37}) leads to the bivariate VGF of queue content and number with the departing batch as
\bea \hspace*{-0.3cm}\boldsymbol{\Pi}^+(z,y)&=& \left[\sum_{n=0}^{a-1}\boldsymbol{\psi}^+(n)\left\{\sum_{j=a}^{b}\boldsymbol{C}_{jn}\textbf{A}^{(j)}(z)\textbf{A}^{(b)}(z)\left(y^b-y^j\right)+\sum_{j=a}^{b}
\boldsymbol{C}_{jn} y^j z^b \textbf{A}^{(j)}(z)\right.\right.\nn\\
&&\left. \left.-  \left(y^b z^n \textbf{I} -y^b\sum_{i=b-a+2}^{b+1-n}\tau_i(z)~\boldsymbol{Q}_{b+1-n-i}~z^{b+1-i}\right) \textbf{A}^{(b)}(z) \right\}\right.\nn\\
&&\left.+\sum_{n=a}^{b-1}\boldsymbol{\psi}^+(n)\left\{ \left(y^b-y^n\right)\textbf{A}^{(n)}(z)\textbf{A}^{(b)}(z)\right.\right.\nn\\&&\left.\left.+\left(y^n z^b \textbf{A}^{(n)}(z)-y^b z^n \textbf{A}^{(b)}(z)\right)\right\}\right] \left[z^b\textbf{I}-\textbf{A}^{(b)}(z)\right]^{-1}\lb{map50}\eea
\vspace*{0.15cm}
\begin{remark}
The above expression represents the bivariate VGF of queue length and size of the departing batch which is the spectrum of the whole analysis. To the best of authors' knowledge no such result is available so far in the literature. Moreover, the similar remark holds for the expression presented in (\ref{map48}) which is the VGF of \emph{only} queue length distribution at departure epoch.
\end{remark}
\subsection{Determination of unknown vectors}
This subsection presents the determination of unknown probability vectors appearing in the bivariate VGF given in (\ref{map50}). In order to determine those unknown vectors, the procedure given in Singh et al. \cite{singh2013computational,singh2016detailed} has been followed. The numerator of bivariate VGF contains $b$ unknown vectors $\{\boldsymbol{\psi}^+(n)\}_{n=0}^{b-1}$, i.e., in total $mb$ unknowns $\{\psi_i^+(n)\}_{n=0}^{b-1},~1\leq i \leq m $, which has to be perceived first. Without loss of any generality, we evaluate the unknown vectors using (\ref{map48}) instead of using (\ref{map50}) as the unknowns appearing in (\ref{map50}) and (\ref{map48}) are exactly the same.\\
\hspace*{0.3cm}The existing literatures clearly identify that the distributions with L.-S.T. as rational function has a considerable impact in applications, for example see Botta et al. \cite{botta1987characterizations}. In view of this, the distributions with rational L.-S.T. of the form $\widetilde{S}_r(\theta)=\frac{P_r(\theta)}{Q_r(\theta)}$ with degree of $P_r(\theta)$ less or equal to that of $Q_r(\theta)$, has to be considered
here. Even if, we can deal with transcendental L.-S.T. (for deterministic distribution) which is rationalized using Pad$\acute{e}$ approximation. Each element of $\textbf{A}^{(r)}(z),~(a \leq r \leq b)$, is also a rational function possessing same denominator, say $d^{(r)}(z)$, because of $\widetilde{S}_r(\theta)$ being a rational function. Consequently, each element of the matrix $\left[z^b\textbf{I}-\textbf{A}^{(b)}(z) \right]$ must be a rational function with the same denominator $d^{(b)}(z)$.\\
\hspace*{0.3cm}Now, we assume that the $(i,j)$-th element of $\textbf{A}^{(r)}(z)$ is $\frac{f^{(r)}_{i,j}(z)}{d^{(r)}(z)},~~1\leq i,j \leq m$. Hence, the $(i,j)$-th element of the matrix $\left[z^b\textbf{I}-\textbf{A}^{(b)}(z) \right]$ is given as
\bea \left[z^b\textbf{I}-\textbf{A}^{(b)}(z) \right]_{i,j}=\frac{h_{i,j}(z)}{d^{(b)}(z)},~~\mbox{where} ~~
h_{i,j} (z)= \left\{\begin{array}{r@{\mskip\thickmuskip}l}
& z^b d^{(b)}(z)-f^{(b)}_{i,i}(z),~~i=j\\
&-f^{(b)}_{i,j}(z),~~~~~~~~~~~~~~i\neq j
 \end{array}\right.\lb{map51}\eea
\\
As both sides of (\ref{map48}) presents row vectors of dimension $(1\times m)$, then comparing element-wise we obtain $m$ simultaneous equations in $m$ unknowns $\Psi_j^+(z),~~1\leq j\leq m$. We get a simplified form of the equations as
\bea H(z) \left[\Psi^+(z)\right]^T= \left[\boldsymbol{ \overline{\chi}}(z)\right]^T\nn\eea
where  $H(z)$ is the square matrices with $(k,l)$-th elements as $\left[H(z)\right]_{k,l}=h_{l,k}(z)$ and
$\boldsymbol{\overline{\chi}}(z)=\left[\boldsymbol{\chi}_1(z),\boldsymbol{\chi}_2(z),\ldots,\boldsymbol{\chi}_m(z)\right]$
and for $1\leq j \leq m$,
\begin{small}
\bea \boldsymbol{\chi}_j(z)&=&\left[\sum_{l=1}^{m} \sum_{n=0}^{a-1}\sum_{k=1}^{m}\psi^+_{k}(n)\left\{ \sum_{\xi=a}^{b}\left(\boldsymbol{C}_{\xi n}\right)_{k,l}\prod_{t=a,~ t\neq \xi}^{b}d^{(t)}(z) f_{l,j}^{(\xi)}(z)\right\}  z^b \right.\nn\\&&\left.+\sum_{l=1}^{m} \sum_{n=0}^{a-1}\sum_{k=1}^{m}\psi^+_{k}(n) \sum_{i=b-a+2}^{b+1-n}\left(\sum_{\gamma=1}^{m}\left\{\tau_i(z)_{k,\gamma} (\boldsymbol{Q}_{b+1-n-i})_{\gamma,l}\right\}z^{b+1-i}\right) f_{l,j}^{(b)}(z) \prod_{t=a}^{b-1}d^{(t)}(z)\right.\nn\\&&\left.-\sum_{l=1}^{m}\sum_{n=0}^{a-1}\psi^+_{l}(n)z^n f_{l,j}^{(b)}(z)\prod_{t=a}^{b-1}d^{(t)}(z) \right.\nn\\
   &&\left.+\sum_{l=1}^{m}\sum_{n=a}^{b-1}\psi^+_{l}(n)\left( z^b d^{(b)}(z) f_{l,j}^{(n)}(z)-z^nd^{(n)}(z) f_{l,j}^{(b)}(z)\right)\prod_{t=a,~ t\neq n}^{b-1}d^{(t)}(z)\right]\nn\\
   &&  \bigg/ \left[~\prod_{t=a}^{b-1}d^{(t)}(z)\right] \lb{map52}\eea
   \end{small}
\\
Using the classical Cramer's rule, the above system of equations can be solved easily. Hence, we are led to $\Psi_j^+(z),~1\leq j \leq m$,
\bea \Psi_j^+(z)=\frac{\left|H_j(z)\right|}{\left|H(z)\right|},~~1\leq j \leq m, \lb{map53}\eea
where both $H_j(z)$ and $H(z)$ are square matrices with $(k,l)$-th elements given by
\bea \left[H_j(z)\right]_{k,l} = \left\{\begin{array}{r@{\mskip\thickmuskip}l}
& h_{l,k}(z),~~~~l\neq j\\
&\boldsymbol{\chi}_{k}(z),~~~~l= j
 \end{array}\right. ~~~~~~\mbox{and}~~~~~~ \left[H(z)\right]_{k,l}=h_{l,k}(z)\lb{map54}\eea
The $j$-th column of the square matrix $H_j(z)$ is replaced by $[\boldsymbol{\chi}_1(z),\boldsymbol{\chi}_2(z),\ldots,\boldsymbol{\chi}_m(z)]^T$ and all other elements are same as those of $H(z)$.\\
Let us assume that $|H(z)|$, which is a polynomial in $z$ must possess a non-zero coefficient of power of $z$. Finally, we have
\bea \Psi_j^+(z)=\frac{\Upsilon_j(z)}{\Upsilon(z)},~~1\leq j \leq m, \lb{map55}\eea
where $\Upsilon_j(z)=|H_j(z)|$ and $\Upsilon(z)=|H(z)|$. More precisely it may be noted that we are having the pgf of \emph{only} queue length distribution for each phase at departure epoch. Now we concentrate on the determination of unknown probability vectors. Consequently, we consider (\ref{map55}), and let us call $\Upsilon(z) = 0$ as characteristic equation associated with the pgf of each phase. It can be proved that $|z^b\textbf{I}-\textbf{A}^{(b)}(z)|\equiv \frac{\Upsilon(z)}{\{d^{(b)}(z)\}^m}=0$ has exactly `$mb$' roots inside and on the closed complex unit disk $|z|\leq 1$, see Gail et al. \cite{gail1995linear}. In this context, we assume that these roots are distinct and denote them as $z_1,z_2,\ldots,z_{mb}$ with $z_{mb}=1$. However, in case of multiple roots the procedure needs a slight modification.\\
\hspace*{0.3cm}The analytical nature of $\Psi_j^+(z)$ in $|z|\leq 1$ implies that the roots $z_1,z_2,\ldots,z_{mb-1}$ of $\Upsilon(z)=0$ (the denominator of (\ref{map55})) must coincide with that of numerator. Now, considering any one component of $\boldsymbol{\Psi}^+(z)$, say $\Psi_j^+(z),~~(1\leq j \leq m)$, we have $mb-1$ equations as
\bea \Upsilon_j(z_i)=0,~~~1\leq i \leq mb-1. \lb{map56}\eea
By employing the normalizing condition $\boldsymbol{\Psi}^+(1)e=1$, we are led to one more equation as
\bea \sum_{j=1}^{m}\Upsilon_j^{'}(1)= \Upsilon^{'}(1)\lb{map57}\eea
Solving (\ref{map56}) and (\ref{map57}) together we obtain `$mb$' unknowns $\psi^+_{j}(n),~~(0\leq n\leq b-1,~1\leq j \leq m)$.
\subsection{Extraction of probability vectors from the VGF}
After determining the unknown vectors $\left\{\boldsymbol{\psi}^+(n)\right\}_{n=0}^{b-1}$~, we change our focus to extract probability vectors $\bpi^+(n,r),~n\geq 0,~a\leq r \leq b$, from completely known bivariate VGF. This can be perceived by inverting $\Pi^+(z,y)$, which is not easily tractable. For this purpose, the coefficient of $y^j,~a\leq j\leq b$, have been accumulated from both the sides of (\ref{map50}) and are precisely given by
\begin{small}
\bea \mbox{coefficient of} ~ y^a:&&~~~~\nn\\
 &&\hspace*{-1.0cm}\sum_{n=0}^{\infty}\bpi^+(n,a)z^n = \left(\sum_{i=0}^{a-1}\boldsymbol{\psi}^+(i)\boldsymbol{C}_{a i}+\boldsymbol{\psi}^+(a)\right)\textbf{A}^{(a)}(z)\lb{map58}\\
\mbox{coefficient of} ~ y^j:&&~~~~\nn\\
&&\hspace*{-1.0cm}\sum_{n=0}^{\infty}\bpi^+(n,j)z^n = \left(\sum_{i=0}^{a-1}\boldsymbol{\psi}^+(i)\boldsymbol{C}_{j i}+\boldsymbol{\psi}^+(j)\right)\textbf{A}^{(j)}(z),~~a+1\leq j \leq b-1\lb{map59}\eea
\bea \hspace*{-1.8cm}\mbox{coefficient of} ~ y^b:~~~~\nn\\
\sum_{n=0}^{\infty}\bpi^+(n,b)z^n &=& \left[\sum_{j=0}^{a-1}\boldsymbol{\psi}^+(j)\left\{\sum_{\ell=a}^{b}\boldsymbol{C}_{\ell j} \textbf{A}^{(\ell)}(z)-\boldsymbol{C}_{b j}  \textbf{A}^{(b)}(z) \right.\right.\nn\\
&&\left.\left. + \left( \boldsymbol{C}_{b j} z^b+\sum_{i=b-a+2}^{b+1-j}\tau_i(z)~ \boldsymbol{Q}_{b+1-j-i}~ z^{b+1-i}-z^j \textbf{I} \right)   \right.\right\}\nn\\
&&\left.+\sum_{j=a}^{b-1}\boldsymbol{\psi}^+(j)\left\{\textbf{A}^{(j)}(z)-z^j \textbf{I}\right\}\right]\textbf{A}^{(b)}(z)
\left[z^b\textbf{I}-\textbf{A}^{(b)}(z)\right]
^{-1}\lb{map60}\eea
\end{small}
Now collecting the coefficient of $z^n$ from both the sides of (\ref{map58}) and (\ref{map59}) we get
\bea \bpi^+(n,a) &=& \left(\sum_{i=0}^{a-1}\boldsymbol{\psi}^+(i)\boldsymbol{C}_{a i}+\boldsymbol{\psi}^+(a)\right)\textbf{A}^{(a)}_n,~~n\geq 0\lb{map61}\\
\bpi^+(n,j) &=& \left(\sum_{i=0}^{a-1}\boldsymbol{\psi}^+(i)\boldsymbol{C}_{j i}+\boldsymbol{\psi}^+(j)\right)\textbf{A}^{(j)}_n,~~a+1\leq j \leq b-1,~~n\geq 0\lb{map62}\eea
It may be noted here that collection of the coefficient of $z^n$ from both the sides of (\ref{map60}) is not an easy task. In order to extract $\bpi^+(n,b)$, we invert (\ref{map60}), where each component of the vector is simply a polynomial in $z$. Let us denote $\sum_{n=0}^{\infty}\bpi^+(n,b)z^n$ as $\boldsymbol{F}^+(z)=[F_1^+(z),\ldots,F_m^+(z) ]$ for simplicity which will be used in rest of the analysis of this section. For the extraction of the probability vectors from $\boldsymbol{F}^+(z)$ the same procedure carried out in the previous section for $\boldsymbol{\Psi}^+(z)$, has to be followed here. In view of this, $\boldsymbol{\Psi}^+(z)$ and $\boldsymbol{\chi}_j(z)$ (used in earlier case in eqn. (\ref{map52})) has to be replaced by $\boldsymbol{F}^+(z)$ and $\Phi_j(z)$, respectively, where $\Phi_j(z)$ is given by\\
\begin{small}
\bea \Phi_{j}(z)&=&\left[\sum_{i=1}^{m}\left\{\sum_{l=1}^{m}\sum_{n=0}^{a-1}\psi^+_{l}(n)\left(\sum_{\xi=a}^{b-1}\boldsymbol{C}_{\xi n}f^{(\xi)}(z) \prod_{t=a, ~t\neq \xi}^{b-1}d^{(t)}(z) \right)_{l,i}\right.\right.\nn\\
&&\left.\left.+\sum_{l=1}^{m}\sum_{n=0}^{a-1}\psi^+_{l}(n)\left(\boldsymbol{C}_{b n} z^b \prod_{t=a}^{b-1}d^{(t)}(z) \right)_{l,i}\right.\right.\nn\\&& \left.\left.+\sum_{l=1}^{m}\sum_{n=0}^{a-1}\psi^+_{l}(n)\left(\sum_{\zeta=b-a+2}^{b+1-n}\tau_{\zeta}(z)\boldsymbol{Q}_{b+1-n-\zeta}~z^{b+1-\zeta} \right)_{l,i}~ \prod_{t=a}^{b-1}d^{(t)}(z)\right.\right.\nn\\
 &&\left.\left.+\sum_{l=1}^{m}\sum_{n=a}^{b-1}\psi^+_{l}(n)f^{(n)}_{l,i}(z)\prod_{t=a,~t\neq n}^{b-1}d^{(t)}(z)\right.\right.\nn\\
&&\left.\left.-\left(\sum_{n=0}^{b-1}\psi^+_{i}(n)z^n\right)\prod_{t=a}^{b-1}d^{(t)}(z)\right\}f^{(b)}_{i,j}(z)\right]
\bigg/\left[\prod_{t=a}^{b-1}d^{(t)}(z)\right],~~~~1\leq j \leq m \lb{map63}\eea
\end{small}
Hence, the simplified form of $F_j^+(z)$ is given by
\bea F_j^+(z)=\frac{\left|G_j(z)\right|}{\left|G(z)\right|},~~1\leq j \leq m, \lb{map64}\eea
where both $G_j(z)$ and $G(z)$ represent square matrix with $(k,l)$-th elements given by
\bea \left[G_j(z)\right]_{k,l} = \left\{\begin{array}{r@{\mskip\thickmuskip}l}
& h_{l,k}(z),~~~~l\neq j\\
&\Phi_{k}(z),~~~~l= j
 \end{array}\right. ~~~~~~\mbox{and}~~~~~~ \left[G(z)\right]_{k,l}=h_{l,k}(z)\lb{map65}\eea
The $j$-th column of the square matrix $G_j(z)$ is replaced by $[\Phi_1(z),\Phi_2(z),\ldots,\Phi_m(z)]^T$ and all other elements are same as those of $G(z)$.\\
We assume that $|G(z)|$, which is a polynomial in $z$, possess a non-zero coefficient of power of $z$. Finally, we have
\bea F_j^+(z)=\frac{\Omega_j(z)}{\Omega(z)},~~1\leq j \leq m, \lb{map66}\eea
where $\Omega_j(z)=|G_j(z)|$ and $\Omega(z)=|G(z)|$.
Now $F_j^+(z)$ being a rational function with completely known polynomials, we can proceed to find it's partial fraction.
Let $\Omega_j(z)$ and $\Omega(z)$ are the polynomials of degree $L_1$ and $M_1$, respectively. Depending upon the distinct and multiple roots of $\Omega(z)$, we discuss possible cases below.
\subsubsection{When all the zeros of $\Omega(z)$ in $|z|>1$ are distinct}
As $\Omega(z)$ possess `$mb$' simple zeros inside and on the unit circle, it is evident that $\Omega(z)$ has $M_1-mb$ distinct zeros in $|z|>1$. Let us denote these zeros as $\gamma_1,\gamma_2,\ldots,\gamma_{M_1-mb}$.\\
\\
\textbf{Case 1:} $L_1\geq M_1$\\
\hspace*{0.5cm}Applying the partial-fraction expansion, we can uniquely write the rational function $F_j^+(z)~(1\leq j \leq m)$ as
\bea F_j^+(z)=\sum_{i=0}^{L_1-M_1}\epsilon_{i,j}z^i +\sum_{k=1}^{M_1-mb}\frac{\eta_{k,j}}{\gamma_{k}-z},\lb{map67}\eea for some constants $\epsilon_{i,j}$ and $\eta_{k,j}$'s. The constants $\epsilon_{i,j}$ in the first summation term are obtained by the division of the polynomial $\Omega_j(z)$ by $\Omega(z)$. Use of the classical residue theorem also leads to
\bee \eta_{k,j}=-\frac{\Omega_j(\gamma_{k})}{\Omega^{\prime}(\gamma_{k})},\quad k=1,2,\ldots,M_1-mb. \eee
Now, collecting the coefficient of $z^n$ from both the sides of (\ref{map67}), we have
\bea \pi_j^+(n,b)&=&\epsilon_{n,j}+\sum\limits_{k=1}^{M_1-mb}\frac{\eta_{k,j}}{\gamma^{n+1}_{k}},\quad n\geq 0. \lb{map68}\eea
\textbf{Case 2:} $L_1< M_1$\\
Use of partial-fraction technique on $F_j^+(z)$ leads to
\bea F_j^+(z)=\sum_{k=1}^{M_1-mb}\frac{\eta_{k,j}}{\gamma_{k}-z},\lb{map69}\eea where
\bee \eta_{k,j}=-\frac{\Omega_j(\gamma_{k})}{\Omega^{\prime}(\gamma_{k})},\quad k=1,2,\ldots,M_1-mb. \eee
Now, collecting the coefficient of $z^n$ from both the sides of (\ref{map69}), we get
\bea \pi_j^+(n,b)&=&\sum\limits_{k=1}^{M_1-mb}\frac{\eta_{k,j}}{\gamma^{n+1}_{k}},\quad n\geq 0. \lb{map70}\eea
\subsubsection{ When some zeros of $\Omega(z)$ in $|z|>1$ are repeated}
There is also a possibility that the denominator $\Omega(z)$ of $F_j^+(z)$ may posses some multiple/repeated zeros whose modulus value are greater than one. Assume that, $\Omega(z)$ has total $\ell$ multiple zeros, say $\beta_1,\beta_2,\ldots,\beta_{\ell}$ with multiplicity $\delta_1,\delta_2,\ldots,\delta_{\ell}$, respectively. Further it is clear that $\Omega(z)$ has total $(M_1-mb-\varsigma)$ distinct zeros, where $\varsigma=\sum_{i=1}^{\ell}\delta_i$, say $\gamma_1,\gamma_2,\ldots,\gamma_{M_1-mb-\varsigma}$.\\
\textbf{Case 1:} $L_1\geq M_1$\\
Applying the partial-fraction method, $F_j^+(z)~(1\leq j \leq m)$ can be uniquely written as
 \bea F_j^+(z)=\sum_{i=0}^{L_1-M_1}\epsilon_{i,j} z^i+\sum_{k=1}^{M_1-mb-\varsigma}\frac{\sigma_{k,j} }{\gamma_k-z}+\sum_{\nu=1}^{\ell}\sum_{i=1}^{\delta_\nu}\frac{\eta_{\nu,i,j}}{(\beta_\nu-z)^{\delta_\nu-i+1}} \lb{map71} \eea
 where
 \begin{small}
\bea \sigma_{k,j}&=&-\frac{\Omega_j(\gamma_k)}{\Omega^{\prime}(\gamma_k)},\quad k=1,2,\ldots,M_1-mb-\varsigma\nn, \\
 \eta_{\nu,i,j}&=&\frac{1}{\left(\delta_{\nu}-i\right)!}~ \lim_{z\rightarrow \beta_{\nu}}\frac{d^{\left(\delta_\nu-i\right)}}{dz^{\left(\delta_\nu-i\right)}}\left[ \frac{(\beta_\nu-z)^{\delta_\nu}~\Omega_j(z)}{\Omega(z)}\right]
~~ \nu=1,2,\ldots,\ell, ~~ i=1,2,\ldots,\delta_\nu. \nn \eea
\end{small}
Now collecting the coefficient of $z^n$ from both the sides of (\ref{map71}), we have
\begin{small}
\bea \pi_j^+(n,b)&=&\epsilon_{n,j}+\sum\limits_{k=1}^{M_1-mb-\varsigma}\frac{\sigma_{k,j}}{\gamma^{n+1}_k}+\sum_{\nu=1}^{\ell}\sum_{i=1}^{\delta_\nu}
\binom{\delta_\nu+n-i}{\delta_\nu-i}
\frac{\eta_{\nu, i,j}}{\beta_\nu^{\delta_\nu+n+1-i}},\quad n\geq 0. \lb{map72}\eea
\end{small}
\textbf{Case 2:} $L_1< M_1$\\
 In this case, in partial-fraction, only the first summation term of the right hand side of (\ref{map71}) has been omitted. Now, collecting the coefficients of $z^n$ one can obtain $\pi_j^+(n,b)$ which are given by
\bea \pi_j^+(n,b)&=&\sum\limits_{k=1}^{M_1-mb-\varsigma}\frac{\sigma_{k,j}}{\gamma^{n+1}_k}+\sum_{\nu=1}^{\ell}\sum_{i=1}^{\delta_\nu}\binom{\delta_\nu+n-i}
{\delta_\nu-i}
\frac{\eta_{\nu, i,j}}{\beta_\nu^{\delta_\nu+n+1-i}},\quad n\geq 0. \lb{map73}\eea
This completes the analysis of obtaining the departure epoch probability vectors. Now, a relation between departure and arbitrary epoch probability vectors is to be established in the next section.
\section{Queue length and server content distribution at arbitrary epoch}\lb{sec5}
In order to obtain the system length distribution and several key performance measures of the concerned queueing model, the joint distribution of queue content and server content at arbitrary epoch plays a pivotal role. We establish a correspondence between departure and arbitrary epoch probability vectors in the following theorem.
\begin{theorem}
The state probability vectors \{ $\textbf{p}(n,0), \bpi(n,r)$\} and \{ $\bpi^+(n,r), \boldsymbol{\psi}^+(n)$ \} are connected by
\begin{small}
\bea \textbf{p}(n,0)&=& \left[\frac{1}{E^*}\sum_{j=0}^{n}\boldsymbol{\psi}^+(j)\boldsymbol{M}_{n,j}\right](-\textbf{D}_0)^{-1},~~0\leq n \leq a-1 \lb{map74}\\
\bpi(0,r)&=& \left[\frac{1}{E^*}\left\{\bpi^+(0,r)-\boldsymbol{\psi}^+(r)  \right\}-\sum_{i=0}^{a-1}\textbf{p}(i,0)\textbf{D}_{r-i}\right]\left(\textbf{D}_0\right)^{-1},~~a\leq r \leq b\lb{map75}\\
\bpi(n,r)&=& \left[\frac{1}{E^*}\bpi^+(n,r)-\sum_{i=1}^{n}\bpi(n-i,r)\textbf{D}_i\right]\left(\textbf{D}_0\right)^{-1},~~a\leq r\leq b-1,~~n\geq 1\lb{map77}\\
\bpi(n,b)&=& \left[\frac{1}{E^*}\left\{\bpi^+(n,b)-\boldsymbol{\psi}^+(n+b) \right\}-\sum_{i=1}^{n}\bpi(n-i,b)\textbf{D}_i\right.\nn\\&&\hspace*{4cm}\left.-\sum_{j=a}^{b}\textbf{p}(b-j,0)\textbf{D}_{n+j}\right]\left(\textbf{D}_0\right)^{-1},~~n\geq 1\lb{map78}\eea
\end{small}
\end{theorem}
where
\begin{small}
\bea E^*&=&\omega+\sum_{n=0}^{a-1}\sum_{j=0}^{n}\boldsymbol{\psi}^+(j)\boldsymbol{M}_{n,j}(-\textbf{D}_0)^{-1}\bec ~~\mbox{and}\nn\\
\omega &=&\sum_{n=0}^{a-1}\boldsymbol{\psi}^+(n) \left(\sum_{\ell=a}^{b}\boldsymbol{C}_{\ell n}s_{\ell}
+\sum_{j=n}^{a-1}\boldsymbol{M}_{j n}\sum_{i=b+1-j}^{\infty}\boldsymbol{\overline{D}}_i s_{b}\right) \bec+\sum_{n=a}^{b}s_n \boldsymbol{\psi}^+(n)\bec+s_b\sum_{n=b+1}^{\infty}\boldsymbol{\psi}^+(n)\bec \nn\eea
\end{small}
\\
\noindent\emph{Proof:} From the equation (\ref{map31}) of lemma 4, we get the desired result (\ref{map74}).
Now setting $\theta=0$ in (\ref{map14}) - (\ref{map17}) we get,
\begin{small}
\bea \hspace*{-1cm}\boldsymbol{\pi}(0,r,0)&=&\boldsymbol{\pi}(0,r)\textbf{D}_0+\sum_{i=0}^{a-1}\textbf{p}(i,0)\textbf{D}_{r-i}+\sum_{j=a}^{b}\boldsymbol{\pi}(r,j,0),~~a\leq r \leq b \lb{map79}\\
\hspace*{-1cm}\boldsymbol{\pi}(n,r,0) &=& \boldsymbol{\pi}(n,r)\textbf{D}_0+\sum_{i=1}^{n}\boldsymbol{\pi}(n-i,r)\textbf{D}_i,~ a\leq r \leq b-1,~n\geq 1 \lb{map80}\\
\hspace*{-1cm}\boldsymbol{\pi}(n,b,0)&=&\boldsymbol{\pi}(n,b)\textbf{D}_0+\sum_{i=1}^{n}\boldsymbol{\pi}(n-i,b)\textbf{D}_{i}+\sum_{j=a}^{b}\textbf{p}(b-j,0)\textbf{D}_{n+j}
+\sum_{j=a}^{b}\boldsymbol{\pi}(n+b,j,0) ,~n\geq 1. \lb{map81} \eea
\end{small}
Dividing (\ref{map79}) - (\ref{map81}) by $\displaystyle \sum_{\ell=0}^{\infty}\sum_{r=a}^{b}\boldsymbol{\pi}(\ell,r,0)e$ and using (\ref{map22}), (\ref{map24}) and lemma 3, after some algebraic simplification we obtain the desired results of (\ref{map75}) - (\ref{map78}).
\section{Queue length and server content distribution at pre-arrival epoch}\lb{sec6}
Having found arbitrary epoch probability vectors, we find an association between arbitrary and pre-arrival epoch probability vectors.
Let $\textbf{p}^-(n,0),~(0\leq n \leq a-1)$ and $\bpi^-(n,r),~(a\leq r \leq b,~n\geq 0)$ be the $1\times m$ vectors with $i$-th  component as $p_i^-(n,0)$ and $\pi_i^-(n,r)$, respectively. Let us define $p_i^-(n,0)$ as the steady-state probability that an arrival finds $n$ $(0\leq n \leq a-1)$ customers in the queue, server idle, and phase of the arrival process is $i$. Similarly, we define $\pi^-(n,r)$ to be the steady-state probability that an arrival finds $n$ $(n\geq 0)$ customers in the queue, server busy with $r$ $(a\leq r \leq b)$ customers and phase of the arrival process is $i$. Then the vectors $\textbf{p}^-(n,0)$ and $\bpi^-(n,r)$ are given by
\bea \displaystyle \textbf{p}^-(n,0)&=&\displaystyle \frac{\textbf{p}(n,0)\displaystyle\sum_{i=1}^{\infty}\textbf{D}_i}{\la_g},~~0\leq n \leq a-1\lb{map83}\\
\bpi^-(n,r)&=&\frac{\bpi(n,r)\displaystyle\sum_{i=1}^{\infty}\textbf{D}_i}{\la_g},~~a\leq r \leq b,~n\geq 0. \lb{map84} \eea
\section{System length distribution and performance measures}\lb{sec7}
Having found the probability vectors $\textbf{p}(n,0)$, ($0\leq n \leq a-1$), $\bpi(n,r)$, ($a\leq r \leq b$, $n\geq 0$), the other significant distribution of interest can be easily obtained and are given below.
\begin{itemize}
 \item Distribution of the number of customers in the system at an arbitrary epoch (including number of customers with the server) is given by\\
$
p_n^{system} = \left\{\begin{array}{r@{\mskip\thickmuskip}l}
& \textbf{p}(n,0)\bec\hspace{3.5cm} 0\leq n \leq a-1,\\
& \displaystyle \sum_{r=a}^{min(b,n)}\bpi(n-r,r)\bec\hspace{1.8cm} a\leq n \leq b,\vspace{0.3cm}\\
&\displaystyle \sum_{r=a}^{b}\bpi(n-r,r)\bec\hspace{2.5cm} n\geq b+1.
\end{array}\right.$\\

\item Distribution of the number of customers in the queue at arbitrary epoch is given by\\
$
p_n^{queue} =\left\{\begin{array}{r@{\mskip\thickmuskip}l}
& \textbf{p}(n,0)\bec+\displaystyle \sum_{r=a}^{b}\bpi(n,r)\bec\hspace{1.4cm} 0\leq n \leq a-1,\\
&\displaystyle \sum_{r=a}^{b}\bpi(n,r)\bec\hspace{3.2cm} n\geq a.
\end{array}\right.
$

\item Distribution of the number of customers in service given that server is busy
\bea p_r^{server} = c\sum_{n=0}^{\infty}\bpi(n,r)\bec,~~  a\leq r \leq b \nn\eea
where $c^{-1}=\left[1-\displaystyle\sum_{n=0}^{a-1}\textbf{p}(n,0)\bec\right]=\mbox{probability that the server is busy} ~(P_{busy}).$
\end{itemize}
It is very much essential to study the performance measures of the queueing system as they play a notable role in designing and improving the efficiency of the system. Some performance measures are listed below:
\begin{itemize}
\item average number of customers waiting in the queue ($L_q)=\displaystyle\sum_{n=0}^{\infty}n p_n^{queue}$,
\item mean number of customers in the system ($L)=\displaystyle\sum_{n=0}^{\infty}n p_n^{system}$,
\item average number of customers with the server ($L_s)=\displaystyle\sum_{r=a}^{b}r p_r^{server}$,
\item mean waiting time of a customer in the queue $(W_q)=\displaystyle \frac{L_q}{\la^*}$, as well as in the system $(W)=\displaystyle \frac{L}{\la^*}$.
\item the probability that the server is idle $(P_{idle})=\displaystyle\sum_{n=0}^{a-1}\textbf{p}(n,0)\bec$.
\end{itemize}
\section{Numerical examples}\lb{sec8}
The feasibility and applicability of the methodology discussed above has been illustrated through some numerical examples as it brings out the inner feelings about the concerned queue to the readers. In this section, we present numerical examples by evoking the service time distributions as phase type (PH) and deterministic ($D$).\\
In the first example, the service time has been considered to follow PH-type distribution which has the representation as ($\boldsymbol{\beta},~\textbf{T}$), where $\boldsymbol{\beta}$ is a row vector of order $\nu$, and $\textbf{T}$ is a square matrix of order $\nu$. The joint queue and server content distribution for $BMAP/G_n^{(7,13)}/1$ queue, where $G$ follows PH distribution, at different epochs (departure and arbitrary) has been displayed in Tables \ref{mapbsd1} - \ref{mapbsd2} with the following input parameters. The $BMAP$ is represented by the matrices\\ \begin{small} $\textbf{D}_0=\bordermatrix{ &    &    &  \cr
                       &-0.542410   &  0.003728 & 0.000000  \cr
                          &0.004349   &  -0.022989 & 0.000621  \cr
                          &0.000000   &  0.001243 & -0.269670  \cr}$, $\textbf{D}_1=\bordermatrix{ &    &    &  \cr
                       &0.010252   &   0.000000 & 0.259089  \cr
                          &0.000000   &  0.008698 & 0.000311  \cr
                          &0.129554   &   0.002485 & 0.002175  \cr}$,\\
$\textbf{D}_3=\bordermatrix{ &    &    &  \cr
                       &0.010252   &   0.000000 & 0.259089  \cr
                          &0.000002  &  0.008698 & 0.000310  \cr
                          &0.129553   &  0.002485 &  0.002175 \cr}$
so that $\boldsymbol{\overline \pi}=[0.171902,0.490074,0.338023]$. \end{small}The PH-distribution is taken as
$\boldsymbol{\beta} = \bordermatrix{&           &             \cr
                                     &  0.20    & 0.80     \cr}$ and
$\textbf{T}=\bordermatrix{ &        &        \cr
                          & -\mu_r  & \mu_r  \cr
                          & 0  & -\mu_r   \cr}$
for $a\leq r \leq b$, where $a=7$, $b=13$, $m=3$, $\mu_r=r\mu,~ (7\leq r \leq 13)$, $\mu=0.035$, $\la^*=0.384331$, $\la_g=0.192166$ and $\rho=0.077971$.\\
Although, the similar results at pre-arrival epoch can generated using the relations developed in section \ref{sec6}, we have not appended the table due to the lack of space.\\
\hspace*{0.3cm}The illustration of numerical example includes the deterministic ($D$) distribution possessing the transcendental L.-S.T. which is rationalized using Pad$\acute{e}$ approximation, say Pad$\acute{e}(k,\ell)$ where $k$ and $\ell$ are the parameters with $k<\ell$ but not $k<<\ell$. The parameters corresponding to this example are given by $a=4$, $b=7$ and $m=2$,
               $\textbf{D}_0=\bordermatrix{ &    &      \cr
                          &-6.937500   &  0.937500  \cr
                         &0.062500    & -0.195800   \cr}$,
                         $\textbf{D}_1=\bordermatrix{ &    &      \cr
                          &5.400000   &  0.000000 \cr
                         &0.000000    & 0.119970  \cr}$,
                         $\textbf{D}_5=\bordermatrix{ &    &      \cr
                          &0.600000  &  0.000000   \cr
                         &0.000000    & 0.01333   \cr}$
so that $\boldsymbol{\overline \pi}=[0.062500, 0.937500]$, $\la^*=0.699956$, $\mu_r=\frac{\mu}{r},~ (4\leq r \leq 7)$ with $\mu=7.0$ and $\rho=0.1$. The queue and server content distribution at different epochs are presented in Tables \ref{mapbsd3} - \ref{mapbsd4}.\\
\begin{table}[h!]
$\vspace{0.1cm}$
\caption{Joint distribution of queue and server content and phase of the arrival process at departure epoch for $BMAP/G^{(7,13)}_n/1$ queue, with $G\sim $PH}\lb{mapbsd1}
\begin{tiny}
\begin{tabular}{|c|ccc|ccc|ccc|}\hline
&\multicolumn{3}{c|}{$r=7$}&\multicolumn{3}{c|}{$r=8$}&\multicolumn{3}{c|}{$r=9$}\\\hline
$n$ & $\pi^+_1(n,7)$ & $\pi^+_2(n,7)$ & $\pi^+_3(n,7)$ & $\pi^+_1(n,8)$ & $\pi^+_2(n,8)$ & $\pi^+_3(n,8)$ & $\pi^+_1(n,9)$ & $\pi^+_2(n,9)$ & $\pi^+_3(n,9)$ \\\hline
0&0.060869&0.028174&0.094486&0.026801&0.011099&0.044618&0.036431&0.013914&0.055249\\
1&0.017491&0.002619&0.034488&0.007925&0.001005&0.014319&0.009502&0.001120&0.018374\\
2&0.006289&0.000636&0.009845&0.002521&0.000230&0.004192&0.003107&0.000255&0.004751\\
3&0.019295&0.002798&0.037995&0.008661&0.001070&0.015643&0.010308&0.001185&0.019911\\
4&0.013212&0.001332&0.020692&0.005273&0.000480&0.008768&0.006472&0.000530&0.009898\\
5&0.005593&0.000555&0.010872&0.002275&0.000200&0.004089&0.002486&0.000201&0.004737\\
10&0.003530&0.000333&0.005590&0.001248&0.000107&0.002077&0.001355&0.000104&0.002095\\
20&0.000352&0.000032&0.000569&0.000099&0.000008&0.000165&0.000085&0.000006&0.000134\\
30&0.000036&0.000003&0.000059&0.000008&0.000000&0.000013&0.000005&0.000000&0.000008\\
40&0.000004&0.000000&0.000006&0.000000&0.000000&0.000001&0.000000&0.000000&0.000000\\
45&0.000001&0.000000&0.000000&0.000000&0.000000&0.000000&0.000000&0.000000&0.000000\\
$\geq$50&0.000000&0.000000&0.000000&0.000000&0.000000&0.000000&0.000000&0.000000&0.000000\\\hline
Total&0.161800&0.039891&0.275838&0.067675&0.015321&0.116074&0.083543&0.018415&0.138606\\\hline
\end{tabular}
$\vspace{0.1cm}$
\\
\begin{tabular}{|c|ccc|ccc|ccc|}\hline
&\multicolumn{3}{c|}{$r=10$}&\multicolumn{3}{c|}{$r=11$}&\multicolumn{3}{c|}{$r=12$}\\\hline
$n$ & $\pi^+_1(n,10)$ & $\pi^+_2(n,10)$ & $\pi^+_3(n,10)$ & $\pi^+_1(n,11)$ & $\pi^+_2(n,11)$ & $\pi^+_3(n,11)$ & $\pi^+_1(n,12)$ & $\pi^+_2(n,12)$ & $\pi^+_3(n,12)$\\\hline
0&0.002336&0.000612&0.005679&0.001593&0.000402&0.004320&0.001509&0.000353&0.003564\\
1&0.000925&0.000078&0.001125&0.000677&0.000052&0.000731&0.000543&0.000039&0.000658\\
2&0.000188&0.000015&0.000435&0.000119&0.000009&0.000302&0.000102&0.000007&0.000231\\
3&0.000995&0.000083&0.001214&0.000724&0.000054&0.000783&0.000577&0.000042&0.000701\\
4&0.000390&0.000031&0.000903&0.000245&0.000018&0.000625&0.000211&0.000014&0.000476\\
5&0.000215&0.000015&0.000272&0.000144&0.000009&0.000163&0.000106&0.000006&0.000133\\
10&0.000075&0.000006&0.000166&0.000043&0.000003&0.000103&0.000033&0.000002&0.000071\\
20&0.000004&0.000000&0.000008&0.000002&0.000000&0.000004&0.000001&0.000000&0.000002\\
30&0.000000&0.000000&0.000000&0.000000&0.000000&0.000000&0.000000&0.000000&0.000000\\
40&0.000000&0.000000&0.000000&0.000000&0.000000&0.000000&0.000000&0.000000&0.000000\\
45&0.000000&0.000000&0.000000&0.000000&0.000000&0.000000&0.000000&0.000000&0.000000\\
$\geq$50&0.000000&0.000000&0.000000&0.000000&0.000000&0.000000&0.000000&0.000000&0.000000\\\hline
Total&0.006062&0.000909&0.011416&0.004122&0.000588&0.008028&0.003517&0.000493&0.006578\\\hline
\end{tabular}
$\vspace{0.1cm}$
\\
\begin{tabular}{|c|ccc|c|}\hline
&\multicolumn{3}{c|}{$r=13$}&\\\hline
$n$ & $\pi^+_1(n,13)$ & $\pi^+_2(n,13)$ & $\pi^+_3(n,13)$ & $\boldsymbol{\psi}_n^+ \bec$\\\hline
0&0.001250&0.000284&0.003290&0.396836\\
1&0.001463&0.000250&0.002805&0.116193\\
2&0.001155&0.000196&0.002540&0.037128\\
3&0.001502&0.000200&0.002501&0.126247\\
4&0.001256&0.000159&0.002441&0.073432\\
5&0.001123&0.000134&0.001901&0.035233\\
10&0.000486&0.000049&0.000906&0.018389\\
20&0.000059&0.000005&0.000107&0.001645\\
30&0.000006&0.000000&0.000010&0.000148\\
40&0.000000&0.000000&0.000001&0.000012\\
45&0.000000&0.000000&0.000000&0.000001\\
$\geq$50&0.000000&0.000000&0.000000&0.000000\\\hline
Total&0.013515&0.001825&0.025776&1.000000\\\hline
\end{tabular}

\end{tiny}
\end{table}

\newpage
\begin{table}[h!]
$\vspace{0.1cm}$
\centering
\caption{Joint distribution of queue and server content and phase of the arrival process at arbitrary epoch for $MAP/G^{(7,13)}_n/1$ queue, with $G\sim $PH}\lb{mapbsd2}
\begin{tiny}
\begin{tabular}{|c|ccc|ccc|ccc|}\hline
&\multicolumn{3}{c|}{$r=0$}&\multicolumn{3}{c|}{$r=7$}&\multicolumn{3}{c|}{$r=8$}\\\hline
$n$ & $p_1(n,0)$ & $p_2(n,0)$&$p_3(n,0)$&$\pi_1(n,7)$ & $\pi_2(n,7)$ & $\pi_3(n,7)$  & $\pi_1(n,8)$ & $\pi_2(n,8)$  & $\pi_3(n,8)$ \\\hline
0&0.012402&0.117421&0.037496&0.014546&0.006573&0.022425&0.005597&0.002265&0.009253\\
1&0.013069&0.062645&0.025276&0.004103&0.000598&0.008087&0.001623&0.000201&0.002931\\
2&0.007717&0.031385&0.016834&0.001459&0.000144&0.002273&0.000510&0.000045&0.000844\\
3&0.017669&0.078114&0.033999&0.004516&0.000639&0.008889&0.001770&0.000213&0.003195\\
4&0.017670&0.069894&0.038054&0.003062&0.000302&0.004773&0.001066&0.000095&0.001765\\
5&0.015042&0.050578&0.028975&0.001279&0.000124&0.002485&0.000454&0.000039&0.000815\\
6&0.017748&0.063873&0.037374&0.002049&0.000199&0.003201&0.000697&0.000061&0.001154\\
10&&&&0.000797&0.000073&0.001258&0.000245&0.000020&0.000408\\
20&&&&0.000077&0.000007&0.000124&0.000019&0.000001&0.000031\\
30&&&&0.000007&0.000000&0.000012&0.000001&0.000000&0.000002\\
40&&&&0.000000&0.000000&0.000001&0.000000&0.000000&0.000000\\
$\geq 45$&&&&0.000000&0.000000&0.000000&0.000000&0.000000&0.000000\\\hline
Total&0.101320&0.473912&0.218010&0.037797&0.009221&0.064151&0.013829&0.003099&0.023622\\\hline
\end{tabular}
$\vspace{0.1cm}$
\\
\begin{tabular}{|c|ccc|ccc|ccc|}\hline
&\multicolumn{3}{c|}{$r=9$}&\multicolumn{3}{c|}{$r=10$}&\multicolumn{3}{c|}{$r=11$}\\\hline
$n$ & $\pi_1(n,9)$&$\pi_1(n,9)$&$\pi_1(n,9)$&$\pi_1(n,10)$ & $\pi_2(n,10)$ & $\pi_3(n,10)$  & $\pi_1(n,11)$ & $\pi_2(n,11)$  & $\pi_3(n,11)$ \\\hline
0&0.006755&0.002523&0.010173&0.000389&0.000099&0.000940&0.000241&0.000059&0.000649\\
1&0.001726&0.000199&0.003337&0.000151&0.000012&0.000184&0.000100&0.000007&0.000108\\
2&0.000557&0.000045&0.000849&0.000030&0.000002&0.000069&0.000017&0.000001&0.000044\\
3&0.001869&0.000210&0.003609&0.000162&0.000013&0.000197&0.000107&0.000008&0.000116\\
4&0.001161&0.000093&0.001766&0.000063&0.000005&0.000145&0.000036&0.000003&0.000091\\
5&0.000439&0.000035&0.000837&0.000034&0.000002&0.000043&0.000020&0.000001&0.000023\\
10&0.000236&0.000018&0.000364&0.000012&0.000000&0.000026&0.000006&0.000000&0.000014\\
20&0.000014&0.000001&0.000022&0.000000&0.000000&0.000001&0.000000&0.000000&0.000000\\
30&0.000000&0.000000&0.000001&0.000000&0.000000&0.000000&0.000000&0.000000&0.000000\\
40&0.000000&0.000000&0.000000&0.000000&0.000000&0.000000&0.000000&0.000000&0.000000\\
$\geq$45&0.000000&0.000000&0.000000&0.000000&0.000000&0.000000&0.000000&0.000000&0.000000\\\hline
Total&0.015180&0.003314&0.025064&0.000989&0.000146&0.001860&0.000611&0.000086&0.001189\\\hline
\end{tabular}
$\vspace{0.1cm}$
\\
\begin{tabular}{|c|ccc|ccc|c|}\hline
&\multicolumn{3}{c|}{$r=12$}&\multicolumn{3}{c|}{$r=13$}&\\\hline
$n$  & $\pi_1(n,12)$ & $\pi_2(n,12)$ &$\pi_3(n,12)$& $\pi_1(n,13)$ & $\pi_2(n,13)$ &$\pi_3(n,13)$&$p^{queue}_n$\\\hline
0&0.000209&0.000048&0.000491&0.000160&0.000036&0.000418&0.251175\\
1&0.000073&0.000005&0.000089&0.000185&0.000031&0.000355&0.125101\\
2&0.000013&0.000000&0.000030&0.000146&0.000024&0.000321&0.063367\\
3&0.000078&0.000005&0.000095&0.000189&0.000024&0.000315&0.156007\\
4&0.000028&0.000001&0.000063&0.000157&0.000019&0.000306&0.140623\\
5&0.000014&0.000000&0.000017&0.000141&0.000017&0.000238&0.101659\\
10&0.000004&0.000000&0.000009&0.000060&0.000006&0.000112&0.003675\\
20&0.000000&0.000000&0.000000&0.000007&0.000000&0.000013&0.000322\\
30&0.000000&0.000000&0.000000&0.000000&0.000000&0.000001&0.000030\\
40&0.000000&0.000000&0.000000&0.000000&0.000000&0.000000&0.000002\\
$\geq$45&0.000000&0.000000&0.000000&0.000000&0.000000&0.000000&0.000000\\\hline
Total&0.000478&0.000066&0.000893&0.001695&0.000227&0.003231&1.000000\\\hline
\multicolumn{8}{|c|}{$L$=4.549595,~~ $L_q$=2.919973,~~ $L_{s}$=7.881826}\\
\multicolumn{8}{|c|}{}\\
\multicolumn{8}{|c|}{$P_{idle}$=0.793243,~~ $W$=11.837688,~~ $W_q$=7.597540}\\\hline
\end{tabular}
\end{tiny}
\end{table}

\newpage
\begin{table}[h!]
$\vspace{0.1cm}$
\caption{Joint distribution of queue and server content and phase of the arrival process at departure epoch for $BMAP/G^{(4,7)}_n/1$ queue with $G\sim D$}
\lb{mapbsd3}
\begin{tiny}
\begin{tabular}{|c|cc|cc|cc|cc|c|}\hline
&\multicolumn{2}{c|}{$r=4$}&\multicolumn{2}{c|}{$r=5$}&\multicolumn{2}{c|}{$r=6$}&\multicolumn{2}{c|}{$r=7$}&\\\hline
$n$ & $\pi^+_1(n,4)$ & $\pi^+_2(n,4)$ & $\pi^+_1(n,5)$ & $\pi^+_2(n,5)$ & $\pi^+_1(n,6)$ & $\pi^+_2(n,6)$  & $\pi^+_1(n,7)$ & $\pi^+_2(n,7)$  &$\boldsymbol{\psi}_n^+\bec$\\\hline
0&0.007936&0.210652&0.000688&0.047289&0.000508&0.042904&0.000392&0.037206&0.347577\\
1&0.020340&0.042719&0.001376&0.007471&0.000899&0.008462&0.000959&0.041759&0.123985\\
2&0.029995&0.020640&0.002197&0.002861&0.001513&0.003629&0.001700&0.023738&0.086275\\
3&0.030352&0.011661&0.002647&0.001759&0.002090&0.002412&0.002794&0.018858&0.072574\\
4&0.023236&0.005903&0.002481&0.001031&0.002309&0.001567&0.003984&0.015934&0.056445\\
5&0.016539&0.007387&0.002039&0.001376&0.002189&0.001879&0.004929&0.014332&0.050671\\
10&0.005384&0.000964&0.000881&0.000215&0.001221&0.000408&0.005602&0.006950&0.021625\\
20&0.000106&0.000011&0.000029&0.000004&0.000066&0.000011&0.001500&0.001057&0.002785\\
30&0.000000&0.000000&0.000000&0.000000&0.000001&0.000000&0.000209&0.000134&0.000345\\
40&0.000000&0.000000&0.000000&0.000000&0.000000&0.000000&0.000026&0.000017&0.000043\\
50&0.000000&0.000000&0.000000&0.000000&0.000000&0.000000&0.000003&0.000002&0.000005\\
$\geq 55$&0.000000&0.000000&0.000000&0.000000&0.000000&0.000000&0.000000&0.000000&0.000000\\\hline
Total&0.193843&0.316233&0.020501&0.064927&0.021392&0.066081&0.082364&0.234659&1.000000\\\hline
\end{tabular}
\end{tiny}
$\vspace{0.2cm}$
\caption{Joint distribution of queue and server content and phase of the arrival process at arbitrary epoch for $BMAP/G^{(4,7)}_n/1$ queue with $G\sim D$}\lb{mapbsd4}
\begin{tiny}
\begin{tabular}{|c|cc|cc|cc|cc|cc|c|}\hline
&\multicolumn{2}{c|}{$r=0$}&\multicolumn{2}{c|}{$r=4$}&\multicolumn{2}{c|}{$r=5$}&\multicolumn{2}{c|}{$r=6$}&\multicolumn{2}{c|}{$r=7$}&\\
\hline
$n$ & $p_1(n,0)$ & $p_2(n,0)$ &$\pi_1(n,4)$ & $\pi_2(n,4)$ & $\pi_1(n,5)$ & $\pi_2(n,5)$  & $\pi_1(n,6)$ & $\pi_2(n,6)$ & $\pi_1(n,7)$ & $\pi_2(n,7)$ &$p^{queue}_n$ \\\hline
0&0.002376&0.243267&0.006163&0.016297&0.000732&0.004762&0.000870&0.005245&0.000996&0.005365&0.286074\\
1&0.004462&0.239295&0.004419&0.001844&0.000547&0.000413&0.000665&0.000593&0.001816&0.005474&0.259530\\
2&0.006055&0.210505&0.002866&0.000694&0.000385&0.000132&0.000490&0.000221&0.001729&0.002635&0.225712\\
3&0.007130&0.186916&0.001646&0.000307&0.000249&0.000065&0.000342&0.000120&0.001563&0.002018&0.200357\\
4&&&0.000832&0.000124&0.000146&0.000031&0.000222&0.000063&0.001358&0.001629&0.004406\\
5&&&0.000862&0.000251&0.000138&0.000059&0.000207&0.000097&0.001232&0.001407&0.004254\\
10&&&0.000135&0.000019&0.000030&0.000005&0.000057&0.000012&0.000715&0.000572&0.001546\\
20&&&0.000001&0.000000&0.000000&0.000000&0.000000&0.000000&0.000116&0.000076&0.000196\\
30&&&0.000000&0.000000&0.000000&0.000000&0.000000&0.000000&0.000015&0.000009&0.000024\\
40&&&0.000000&0.000000&0.000000&0.000000&0.000000&0.000000&0.000002&0.000001&0.000003\\
50&&&0.000000&0.000000&0.000000&0.000000&0.000000&0.000000&0.000000&0.000000&0.000000\\
$\geq 55$&&&0.000000&0.000000&0.000000&0.000000&0.000000&0.000000&0.000000&0.000000&0.000000\\\hline
Total&0.203023&0.150172&0.300113&0.227053&0.022495&0.016902&0.009456&0.007111&0.021795&0.016398&1.000000\\\hline
\multicolumn{12}{|c|}{$L$=2.108631,~~ $L_q$=1.552581,~~ $L_{s}$=5.560853}\\
\multicolumn{12}{|c|}{}\\
\multicolumn{12}{|c|}{$P_{idle}$=0.900006,~~ $W$=3.012518,~~ $W_q$=2.218111}\\\hline
\end{tabular}
\end{tiny}
\end{table}
\section{Conclusion}\lb{sec9}
In this paper, the supplementary variable technique and bivariate probability generating function approach is adopted to analyze a generally distributed batch-size-dependent service queue with batch Markovian arrival process. Three significant features have been included, firstly, the closed-form expression of bivariate vector generating function of queue and server content distribution at departure epoch, secondly, the extraction of probability vectors, and finally, the relation between probability vectors at departure and arbitrary epochs. We have dealt with several complicated analytic expressions during the analysis. Through assorted numerical examples, it is clear that the methodology is tractable and easily implementable. As a final conclusion, we feel that this queuing model can be used to quantify the effects of multimedia services over a wireless local communication networks.

\bibliographystyle{unsrtnat}
\bibliography{BMAPbib}

\end{document}